\documentclass[12pt]{elsarticle}

\usepackage{mathstyle00}
\usepackage{amsfonts,bm} 
\usepackage{verbatim} 

\usepackage[margin=0.95in]{geometry}    

\begin{document}
\begin{frontmatter}
\title{Sharp error estimates on a stochastic structure-preserving scheme in computing effective diffusivity of 3D chaotic flows}
		
		\author[hku]{Zhongjian Wang}
		\ead{ariswang@connect.hku.hk}
		\author[uci]{Jack Xin}
		\ead{jxin@math.uci.edu}
		\author[hku]{Zhiwen Zhang\corref{cor1}}
		\ead{zhangzw@hku.hk}
		
		\address[hku]{Department of Mathematics, The University of Hong Kong, Pokfulam Road, Hong Kong SAR.}
		\address[uci]{Department of Mathematics, University of California at Irvine, Irvine, CA 92697, USA.}
		\cortext[cor1]{Corresponding author}
\begin{abstract}
	
In this paper, we study the problem of computing the effective diffusivity for particles moving in chaotic flows. Instead of solving a convection-diffusion type cell problem in the Eulerian formulation (arising from homogenization theory for parabolic equations), we compute the motion of particles in the Lagrangian formulation, which is modeled by stochastic differential equations (SDEs). A robust numerical integrator based on a splitting method was proposed to solve the SDEs and a rigorous error analysis for the numerical integrator was provided using the backward error analysis (BEA) technique \cite{WangXinZhang:18}. However, the upper bound in the error estimate is not sharp. To improve our result, we propose a new and uniform in time error analysis for the numerical integrator that allows us to get rid of the exponential growth factor in our previous error estimate. Our new error analysis is based on a probabilistic approach, which interprets the solution process generated by our numerical integrator as a Markov process. By exploring the ergodicity of the solution process, we prove the convergence analysis of our method in computing effective diffusivity over infinite time. We present numerical results to verify the accuracy and efficiency of the proposed method in computing effective diffusivity for several chaotic flows, especially the Arnold-Beltrami-Childress (ABC) flow and Kolmogorov flow in three-dimensional space. 

\end{abstract}
		
\end{frontmatter}

\section{Introduction}\label{sec:introduction}
\noindent
Diffusion enhancement in fluid advection is a fundamental problem to characterize and quantify the large-scale effective diffusion in fluid flows containing complex and turbulent streamlines, which is of great theoretical and practical importance; see e.g. \cite{Fannjiang:94,Fannjiang:97,
	Carmona1997homogenization,mclaughlin1997effect,Yaulandim:1998,Majda:99,
	PavliotisStuart:03,BenOwhadi2003,PavliotisStuart:05,PavliotisStuart:07,
	JackXin:11,JackXin:15}
 and references therein. Its applications can be found in many physical and engineering sciences, including 
atmosphere science, ocean science, chemical engineering, and combustion.  To study the diffusion enhancement phenomenon, one can consider a passive tracer model, which describes particle motion with zero inertia

\begin{align}
	d\textbf{X}(t) = \textbf{v}(\textbf{X},t) + \sigma d\textbf{W}(t),  \quad  \textbf{X}\in \mathbb{R}^{d},   \label{passivemodel}
\end{align}
where $\textbf{X}$ is the position of the particle, $\sigma>0$ is the molecular diffusion coefficient, and $\textbf{W}(t)$ is a $d$-dimensional Brownian motion. The velocity $\textbf{v}(\textbf{X},t)$ satisfies either the Euler or the Navier-Stokes equation. In practice, $\textbf{v}(\textbf{X},t)$ can be modeled by a random field that mimics the energy spectra of the turbulent flow \cite{Majda:99}. 

For spatial-temporal periodic velocity fields and random velocity fields with short-range correlations,
the homogenization theory \cite{BensoussanLionsPapa:2011,Garnier:97,Oleinik:94,Stuart:08} states that the long-time large-scale behavior of the particles is governed by a Brownian motion. More precisely, let $D^{E}\in R^{d\times d}$ denote the effective diffusivity matrix and
$\textbf{X}^{\epsilon}(t)\equiv\epsilon \textbf{X}(t/\epsilon^2)$. Then, $\textbf{X}^{\epsilon}(t)$ converges in distribution to a Brownian motion $\textbf{W}(t)$ with covariance matrix $D^{E}$, i.e., $\textbf{X}^{\epsilon}(t) \xrightarrow{\text{d}}\sqrt{2D^E}\textbf{W}(t)$, as $\epsilon\to 0$. The effective diffusivity matrix $D^E$ can be expressed in terms of particle ensemble average (Lagrangian framework) or integration of solutions to cell problems (Eulerian framework). The dependence of $D^E$ on the velocity field of the problem is highly nontrivial. For {time-independent} Taylor-Green velocity field, the authors of \cite{StuartZygalakis:09} proposed a stochastic splitting method and calculated the effective diffusivity in the limit of vanishing molecular diffusion. For random velocity fields with long-range correlations, various forms of anomalous diffusion, such as super-diffusion and sub-diffusion, can be obtained for exactly solvable models (see \cite{Majda:99} for a review). However, the long-time large-scale behavior of the particle motion is in general difficult to study analytically.

In recent work \cite{WangXinZhang:18}, we proposed a numerical integrator to compute the effective diffusivity 
of chaotic and stochastic flows using structure-preserving schemes. We also investigated the existence of residual diffusivity for several different velocity fields, including the time periodic cellular flows. The residual diffusivity, a special yet remarkable convection-enhanced diffusion phenomenon, refers to the non-zero and finite effective diffusivity in the limit of zero molecular diffusivity as a result of a fully chaotic mixing of the streamlines. Mathematically, we provided a rigorous error estimate for the numerical methods in computing the effective diffusivity. Specifically, let $D^E$ denote the exact effective diffusivity matrix and $D^{E,num}$ denote the numerical result obtained using our method (see the formula in Eq.\eqref{Def_EffectiveDiffusivity_Lagrangian}), respectively. We obtained the error estimate, $|D^{E,num}-D^E|\le C\Delta t+C(T)\Delta t^2$, where the $T$ should be greater than the mixing time. To the best of our knowledge, this result is the first one in the literature to study the convergence on the numerical approximation of the effective diffusivity of chaotic 
flows, which shows that the main source of error does not depend on time. However, the prefactor $C(T)$ in the second term may grow exponentially fast, which makes the estimate not sharp.  

To get a sharp error estimate, we shall develop a new methodology in this paper, which allows us to get rid of the exponential growth factor $C(T)$. Our analysis is based on a probabilistic approach. We interpret the solution process generated by our numerical integrator as a Markov process, where the transition kernel can be constructed explicitly due to the additive noise in the passive tracer model \eqref{passivemodel}. By exploring the ergodicity of the solution process, we succeed in the convergence analysis of our method and give a sharp error estimate for the numerical solution of the effective diffusivity. Most importantly, our convergence analysis reveals the ergodic structure of the solution process, so that {we can compute long-time integration of the passive tracer model in order to accurately compute the effective diffusivity. As we will prove in Theorem \ref{thm:convergence} the error term of the effective diffusivity does not depend on the computational time; see Fig.\ref{fig:eg7f2}}. Finally, we present numerical results to verify the accuracy of the proposed method in computing effective diffusivity for several typical chaotic flow problems of physical interests, including the Arnold-Beltrami-Childress (ABC) flow and the Kolmogorov flow in three-dimensional space. The phenomenon of convection-enhanced diffusion for those velocity fields will also be investigated. 

Our computation of convection-enhanced diffusivity in three-dimensional chaotic flows appears to be the first in the Lagrangian framework. Alternative computation in the Eulerian framework involves singularly perturbed advection-diffusion equations whose solutions develop sharp boundary layers with unknown locations a-priori. We are 
aware of only \cite{Biferale:95} on ABC flows, which we recover and go beyond by two orders of magnitude of molecular diffusivity; see the numerical results in Section \ref{subsec:ExperimentOnPhenomenon} later.

The rest of the paper is organized as follows. In Section 2, we shall review the background of the passive tracer model  and the definition of the effective diffusivity matrix using the Eulerian framework and the Lagrangian framework. 
In Section 3, we propose our numerical integrator in computing the passive tracer model. 
Section 4 is the main part of this paper, where we shall provide our new error estimate based on a probabilistic approach. In addition, we shall show that our method can be used to solve high-dimensional flow problems and the error estimate can be obtained in a straightforward way. In Section 5, we present numerical results to demonstrate the accuracy and efficiency of our method. We also investigate the convection-enhanced diffusivity for several chaotic 
velocity fields, especially the three-dimensional cases. Concluding remarks are made in Section 6.

\section{The definitions of effective diffusivity} \label{sec:EffectiveDiffusivity}
\noindent
We first introduce the definitions of effective diffusivity for chaotic flows. {To be consistent with the setting of the main results in this paper, we assume that the velocity $v$ in Eq.\eqref{passivemodel} is time-independent. Then the SDE \eqref{passivemodel} can be simplified to,}

\begin{align}\label{eqn:generalSDEDefD}
d\textbf{X}(t) = \textbf{v}(\textbf{X}) + \sigma d\textbf{W}(t),  \quad  \textbf{X}\in \mathbb{R}^{d},   
\end{align}
where $\sigma>0$ is the molecular diffusion coefficient, $\textbf{X}$ is the position of the particle, $\textbf{v}(\textbf{X})$ is the Eulerian velocity field at position $\textbf{X}$, $\textbf{W}(t)$ is 
a $d$-dimensional Brownian motion. 
The interested reader is referred to \cite{Biferale:95, Majda:99, PavliotisStuart:03, WangXinZhang:18} and references therein for the results of passive tracer models with time-dependent velocities.   

There are two main frameworks to compute the effective diffusivity of the passive tracer models. We first 
discuss the Eulerian framework. One natural way to study the expectation of the paths for the SDE given by the Eq.\eqref{eqn:generalSDEDefD} is to consider its associated backward Kolmogorov equation. Specifically, given 
a sufficiently smooth function $\phi(\textbf{x})$ in $\mathbb{R}^{d}$, let $u(\textbf{x},t)=\mathbb{E}[\phi(\textbf{X}_t)|\textbf{X}_0=\textbf{x}]$ and $\textbf{X}_t=(x_1(t),...,x_d(t))^{T}$ is the solution to Eq.\eqref{eqn:generalSDEDefD}, then $u(\textbf{x},t)$ satisfies
the backward Kolmogorov equation as
 
\begin{align}\label{BackwardKolmolgorovEquation0}
u_{t}=\mathcal{L}u, \quad u(\textbf{x},0)=\phi(\textbf{x}).
\end{align} 
In Eq.\eqref{BackwardKolmolgorovEquation0}, the generator $\mathcal{L}$ is defined as 

\begin{align}\label{DefinitionGenerator}
\mathcal{L}u = \textbf{v}\cdot \nabla u + D_{0}  \Delta u,
\end{align}
where $D_0=\sigma^2/2$ is the diffusion coefficient and $v$ is the velocity field. When $\textbf{v}(\textbf{x})$ is incompressible (i.e. $\nabla_{\textbf{x}}\cdot \textbf{v}(\textbf{x})=0$), deterministic and periodic in $O(1)$ scale, where we assume the period of $\textbf{v}(\textbf{x})$ is $1$ in each dimension of the physical space, the formula for the effective diffusivity matrix is \cite{BensoussanLionsPapa:2011,Stuart:08}

\begin{align}
D^{E} = D_0I - \big\langle \textbf{v}(\textbf{x})\otimes \chi(\textbf{x}) \big\rangle_{p},
\label{Def_EffectiveDiffusivity_Euler}
\end{align}
where we have assumed that the fluid velocity $\textbf{v}(\textbf{x})$ is smooth and 
the (vector) corrector filed $\chi(\textbf{x})$ satisfies the cell problem,

\begin{align}
- D_0 \Delta  \chi - \textbf{v}(\textbf{x}) \cdot \nabla \chi =  \textbf{v}(\textbf{x}), \quad \textbf{x}\in \mathbb{T}^d,
\label{CellProblem_EffectiveDiffusivity} 
\end{align}
and $\langle \cdot \rangle_{p} $ denotes spatial average over $\mathbb{T}^{d}$. 
Since $v(\textbf{x})$ is incompressible, the solution $\chi(\textbf{x})$ to the cell problem 
\eqref{CellProblem_EffectiveDiffusivity} is unique up to an additive constant by the Fredholm alternative. By multiplying $\chi$ to Eq.\eqref{CellProblem_EffectiveDiffusivity} and integrating in $\mathbb{T}^d$ with consideration of periodicity of $\chi$ and $v$, we will get another equivalent formula for the effective diffusivity,

\begin{align}
D^{E} = D_0I + D_0 \big\langle \nabla\chi(\textbf{x})\otimes \nabla\chi(\textbf{x}) \big\rangle_{p}.
\label{CellProblem_EffectiveDiffusivity2}
\end{align}
The correction to $D_0$ is nonnegative definite in Eq.\eqref{CellProblem_EffectiveDiffusivity2}. We can see that 
$\textbf{e}^{T}D^E\textbf{e}\geq D_0 $ for all unit column vectors $ \textbf{e} \in \mathbb{R}^{d}$, which is called convection-enhanced diffusion. {By energy estimate of $\chi$, one can find an upper bound for the effective diffusivity, i.e., 
for any nonzero unit column vector $\textbf{e}\in \mathbb{R}^{d}$, we have, 
\begin{equation}\label{eqn:maximaldiffusion}
\textbf{e}^{T}D^E\textbf{e} \leq \frac{c}{D_0}, \quad \text{as } D_0\to 0,
\end{equation} 
where the constant $c$ depends on the flow but not on $D_0$. More details of the derivation can be found in \cite{Biferale:95,mezic1996maximal,Fannjiang:94}. We are interested in studying the different scaling laws (between $D_0$ and $\frac{1}{D_0}$) of the convection-enhanced diffusion phenomenon for different chaotic flows in this paper.} The residual diffusivity phenomenon that we studied in \cite{WangXinZhang:18} is one case. While the upper bound given by Eq.\eqref{eqn:maximaldiffusion} is another case, which is called convection-enhanced diffusion with maximal enhancement \cite{mezic1996maximal}; see Fig.\ref{fig:eg7f1} for the result of the ABC flow obtained using our method. 

In practice, the cell problem \eqref{CellProblem_EffectiveDiffusivity} can be solved using numerical methods,
such as spectral methods. In \cite{JackXinLyu:2017}, a small set of adaptive basis functions were constructed from
fully resolved spectral solutions to reduce the computation cost. However, when $D_0$ becomes extremely small,
the solutions of Eq.\eqref{CellProblem_EffectiveDiffusivity} develop sharp gradients and demand a large number of Fourier modes to resolve, which makes the spectral method computationally expensive and unstable.
\begin{remark}
	One can define the adjoint operator $\mathcal{L}^{*}$ as $\mathcal{L}^{*}\rho =-\nabla\cdot(\textbf{v} \rho)+D_{0} \Delta \rho$.  Let $\rho(\textbf{x},t)$ denote the density function of the particle $\textbf{X}(t)$ of Eq.\eqref{eqn:generalSDEDefD}. Then,  $\rho(\textbf{x},t)$	satisfies the Fokker-Planck equation 
	$\rho_{t}=\mathcal{L}^{*} \rho$ with the initial density $\rho(\textbf{x},0)=\rho_0(\textbf{x})$, 
	where $\rho_0(\textbf{x})$ is the density of the particle $\textbf{X}(0)$.
\end{remark}
Alternatively, one can use the Lagrangian framework to compute the effective diffusivity matrix, which is defined by (equivalent to Eq.\eqref{Def_EffectiveDiffusivity_Euler} via the homogenization theory)

\begin{align}
D_{ij}^{E}=\lim_{t\rightarrow\infty}\frac{\Big\langle\big(x_i(t)-x_i(0))(x_j(t)-x_j(0)\big)\Big\rangle}{2t},
\quad 1\leq i,j \leq d,
\label{Def_EffectiveDiffusivity_Lagrangian}
\end{align}
where $\textbf{X}(t)=(x_1(t),...,x_d(t))^{T}$ is the position of a particle tracer at time $t$ and the average $\langle\cdot\rangle $ is taken over an ensemble of test particles.
If the above limit exists, that means the transport of the particle is a standard diffusion process, at least on a long-time scale. If the passive tracer model has a deterministic divergence-free and periodic velocity field, 
this is the typical situation, i.e., the spreading of the particle $\Big\langle\big(x_i(t)-x_i(0))(x_j(t)-x_j(0)\big)\Big\rangle$ grows linearly with respect to the time $t$. For example when the velocity field is given by the Taylor-Green velocity field \cite{Fannjiang:94,StuartZygalakis:09}, the long-time and large-scale behavior of the passive tracer model is a diffusion process. However, there are also cases showing that the spreading of particles does not grow linearly with time but has a power law $t^{\gamma}$, where $\gamma>1$ and $\gamma<1$ correspond to super-diffusive and sub-diffusive behaviors, respectively; see e.g. \cite{Biferale:95,Majda:99,BenOwhadi2003}.

We shall consider the Lagrangian approach in this paper. The Lagrangian framework has the advantages that: 
(1) it is easy to implement; (2) its computational cost linearly depends on the dimension of the passive tracer model; and (3) it does not directly suffer from a small molecular diffusion coefficient $\sigma$ during the computation. 
However, we should point out that the major difficulty in solving Eq.\eqref{eqn:generalSDEDefD} comes from the fact that the computational time should be long enough to approach the diffusion (mixing) time scale. To address this challenge, we shall develop robust numerical integrators, which are structure-preserving and accurate for long-time integration. Moreover, we aim to develop the convergence analysis of the proposed numerical integrators in long-time integration. Finally, we shall investigate the relationship between several typical chaotic flows and the corresponding effective diffusivity.  

\section{Symplectic stochastic integrators}\label{sec:NewStochasticIntegrators}
\subsection{Derivation of numerical integrators}\label{sec:DerivationSchemes}
\noindent
To demonstrate the main idea, we first construct a symplectic stochastic integrator for a two-dimensional passive tracer model with a separable Hamiltonian. { High-dimensional models, including the cases when the velocity field is given by ABC flow and Kolmogorov flow, will be discussed in Section \ref{sec:HighDimensionalCases}.} Specifically, {let $\textbf{X}=(x_1,x_2)^T$ denote the position of the particle and $\textbf{v}=(-f(\textbf{X}),g(\textbf{X}))^T=(-f(x_1,x_2),g(x_1,x_2))^T$ denote the velocity field}, then the passive tracer model can be written as  

\begin{equation}\label{eqn:particleSDE}
\begin{cases}
dx_1=-f(x_1,x_2)dt+\sigma {d} W_{1,t}, \quad  x_1(0)=x_{1}^{0}, \\
dx_2=g(x_1,x_2)dt+\sigma {d}W_{2,t},  \quad \quad   x_2(0)=x_2^{0}, 
\end{cases}
\end{equation}
where $W_{i,t}$, $i=1,2$, are independent Brownian motions. 

{Since the velocity $\textbf{v}$ is generated from a separable Hamiltonian function, we assume that there exists a separable function $H(x_1,x_2)=F(x_2)+G(x_1)$ such that $f(x_1,x_2)=H_{x_2}(x_1,x_2)$, $g(x_1,x_2)=H_{x_1}(x_1,x_2)$, and $H(x_1,x_2)$ is a periodic function on $\mathbb{R}^2$ with period $1$.
We denote with slightly abuse of notation by $f(x_2)$ and $g(x_1)$ for each component of the velocity $\textbf{v}$, i.e., $f(x_2)=f(x_1,x_2)$ and $g(x_1)=g(x_1,x_2)$. These notations simplify our derivation. 
Whenever a statement corresponds to $f(x_2)$ (or $g(x_1)$) is made, it is equivalent to that for $f(x_1,x_2)$ or $g(x_1,x_2)$}. Furthermore, we assume that $H(x_1,x_2)$ is smooth so the first-order derivatives of $f(x_2)$ and $g(x_1)$ are bounded, which guarantee the existence and uniqueness of the solution $(x_1,x_2)$ to the SDE \eqref{eqn:particleSDE}. {The Hamiltonian function is also referred to as the stream function in the fluid mechanical literature.}

In \cite{WangXinZhang:18}, we proposed a structure-preserving scheme based on a Lie-Trotter splitting idea to solve the SDE \eqref{eqn:particleSDE}. Specifically, we split the Eq.\eqref{eqn:particleSDE} into a deterministic subproblem, 

\begin{equation}\label{eqn:particleSDE_1}
\begin{cases}
dx_1=-f(x_2)dt, \\
dx_2=g(x_1)dt, 
\end{cases}
\end{equation} which is solved using a symplectic-preserving scheme (the symplectic Euler scheme for deterministic equations) and a stochastic subproblem,  
\begin{equation}\label{eqn:particleSDE_2}
\begin{cases}
dx_1=\sigma {d}W_{1,t},  \\
dx_2=\sigma {d}W_{2,t}, 
\end{cases}
\end{equation} which is solved using the {Euler-Maruyama scheme} \cite{Oksendal:13}. Eventually, the one step integrator of Eq.\eqref{eqn:particleSDE} is given by, 

\begin{equation}\label{scheme}
\begin{cases}
x_1^{n}=x_1^{n-1}-f(x_2^{n-1})\Delta t+\sigma\sqrt{\Delta t}\xi_1,\\
x_2^{n}=x_2^{n-1}+g\big(x_1^{n-1}-f(x_2^{n-1})\Delta t\big)\Delta t+\sigma\sqrt{\Delta t}\xi_2,
\end{cases}
\end{equation}
where $\xi_1, \xi_2\sim\mathcal{N}(0,1)$ are i.i.d. normal random variables. We denote the stochastic process generated by \eqref{scheme} as $\textbf{X}^n=(x_1^{n},x_2^{n})^T$, which is the numerical approximation to the exact solution {$\textbf{X}(t_n)$ to the SDE \eqref{eqn:particleSDE}  at each lattice point of time $t_n=n\Delta t$}. 

When the Hamiltonian system contains additive temporal noise, the noise itself is considered to be symplectic pathwise \cite{Milstein:02}. We state that the scheme \eqref{scheme} is stochastic symplectic-preserving since it preserves symplecticity as a composition of symplectic transforms and it converges as time-step tends to zero.
Though there are several prior works on developing symplectic-preserving scheme for solving ODEs and PDEs (see \cite{ErnstLubich:06,hong2006multi,JanHesthaven2017structure} and references therein), the novelty of our work is the rigorous theory and sharp estimate on the numerical error in computing the effective diffusivity.

\begin{remark}
In general, the second-order Strang splitting \cite{strang:68} is more frequently adopted to solve ODEs and PDEs. The only difference between the Strang splitting method and the Lie-Trotter splitting method is that the first and last steps are modified by half of the time-step $\Delta t$. For the SDEs, however, the dominant source of error comes from the random subproblem \eqref{eqn:particleSDE_2}. Thus, it is not necessary to implement the Strang splitting scheme here. 
\end{remark}

\begin{remark}
The long-time integration for stochastic Langevin equation was studied in the literature; see e.g.  \cite{bou2010long,abdulle2015long}. However,  passive tracer model \eqref{passivemodel} or \eqref{eqn:particleSDE} studied here has several different features. First, our model problem does not have a damping term so its dynamic behavior and invariant measure of the system are totally different.  In addition, the quantity of interests is different. One of the main focuses in \cite{bou2010long,abdulle2015long} is to investigate whether the average energy remains bounded. 
Our aim here is to study whether the effective diffusivity exists; see the definition in Eq.\eqref{Def_EffectiveDiffusivity_Lagrangian}, and to investigate the convection-enhanced diffusion phenomenon; see Section \ref{subsec:ExperimentOnPhenomenon}.  
\end{remark}

\subsection{The backward Kolmogorov equation and related results}\label{sec:KnownFacts}
\noindent
For the convenience of the reader, we first give a brief review of the theoretical results for the scheme \eqref{scheme} obtained in \cite{WangXinZhang:18} and references therein. We first define the backward Kolmogorov equation associated with the Eq.\eqref{eqn:particleSDE} as 

\begin{align}\label{BackwardKolmolgorovEquation}
u_{t}=\mathcal{L}u, \quad u(\textbf{x},0)=u_0(\textbf{x}),
\end{align}
where the generator $\mathcal{L}$ (associated with the Markov process in Eq. \eqref{eqn:particleSDE}) 
is given by

\begin{align}
\mathcal{L}=-f\partial_{x_1}+g\partial_{x_2}+\frac{1}{2}\sigma^2\partial_{x_1x_1}^2+\frac{1}{2}\sigma^2\partial_{x_2x_2}^2.
\label{HamiltonianFlowOperator}
\end{align}
Recall that the solution $u(\textbf{x},t)$ to the Eq.\eqref{BackwardKolmolgorovEquation} satisfies $u(\textbf{x},t)=\mathbb{E}[\phi(\textbf{X}_t)|\textbf{X}_0=\textbf{x}]$, where $\textbf{X}_t=(x_1(t),x_2(t))^{T}$ is the solution to Eq.\eqref{eqn:particleSDE} and
$\phi$ is a smooth function in $\mathbb{R}^{2}$. 

Similarly, 
 we can study the flow generated by the symplectic splitting scheme \eqref{scheme}.
Recalling the splitting method during the derivation of the scheme in Section \ref{sec:DerivationSchemes}, we define $\mathcal{L}_1=-f\partial_{x_1}$, $\mathcal{L}_2=g\partial_{x_2}$  and $\mathcal{L}_3=\frac{\sigma^2}{2}(\partial^{2}_{x_1x_1}+\partial^{2}_{x_2x_2})$. Starting from $u(\cdot,0)$, we compute

\begin{equation}\label{eqn:flowpde}
\begin{cases}
\partial_tu^1&=\mathcal{L}_1 u^1,\quad u^1(\cdot,0)=u(\cdot,0),\\
\partial_t u^2&=\mathcal{L}_2 u^2,\quad u^2(\cdot,0)=u^1(\cdot,\Delta t), \\
\partial_t u^3&=\mathcal{L}_3u^3,\quad u^3(\cdot,0)=u^2(\cdot,\Delta t).\\
\end{cases}		
\end{equation}
Then $u^3(\cdot,\Delta t)$ will be the flow at time $t=\Delta t$ generated by our scheme and it approximates the solution $u(\cdot,\Delta t)$ to the Eq.\eqref{BackwardKolmolgorovEquation}. It is also worth mentioning that, $u^2(\cdot,\Delta t)$ is the exact flow generated by deterministic symplectic Euler scheme in solving Eq.\eqref{eqn:particleSDE_1}. {And $u^3(\cdot, \Delta t)$ is the flow generated by Euler-Maruyama scheme starting from $u^2(\cdot, \Delta t)$. The latter is due to the fact that Euler-Maruyama schemes are exact when solving white noise SDE like Eq.\eqref{eqn:particleSDE_2}}. Later on, we repeat this process to compute the flow equations of our scheme at 
other time steps, which approximate $u(\cdot,n\Delta t), n=2,3,...$. 

To analyze the error between the flow operator in Eq.\eqref{BackwardKolmolgorovEquation} and {the composition of operators} in Eq.\eqref{eqn:flowpde}, we shall resort to the Baker-Campbell-Hausdorff (BCH) formula, which is widely used in {non-commutative} algebra \cite{BCHformula1974baker}. For example, in the matrix theory,

\begin{equation}\label{eqn:BCHformula}
{\exp(tA )\exp(tB)}=\exp\bigg(t(A+B)+t^2\frac{[A,B]}{2}+\frac{t^3}{12}\Big(\big[A,[A,B]\big]+\big[B,[B,A]\big]\Big)+\cdots\bigg),
\end{equation}
where $t$ is a scalar, $A$ and $B$ are two square matrices with the same size, $[,]$ is the Lie-Bracket,
and the remaining terms on the right hand side are all nested Lie-brackets. 
In our analysis, we replace the matrices in Eq.\eqref{eqn:BCHformula} by differential operators and the BCH formula yields the local structure of our splitting scheme. Let $I_{\Delta t}$ denote the {composite} flow operator associated with Eq.\eqref{eqn:flowpde}, i.e., 

\begin{equation}\label{eqn:operatorflow}
I_{\Delta t} u(\cdot,0):=\exp(\Delta t\mathcal{L}_3)\exp(\Delta t\mathcal{L}_2)\exp(\Delta t\mathcal{L}_1)u(\cdot,0).
\end{equation}
Recall that the exact solution to the Eq.\eqref{BackwardKolmolgorovEquation} at time $t=\Delta t$ can be represented as 

\begin{equation}\label{eqn:TrueFlowOperator}
u(\cdot,\Delta t)=\exp(\Delta t\mathcal{L})u(\cdot,0)=
\exp(\Delta t(\mathcal{L}_{1}+\mathcal{L}_{2}+\mathcal{L}_{3}))u(\cdot,0),
\end{equation}
{or equivalently,
$\mathbb{E}[\textbf{X}_1|\textbf{X}_0=\textbf{x}]=I_{\Delta t}\phi (\textbf{x})$, where expectation are taken over randomness from noise in the scheme \eqref{scheme}.} Now we can apply the BCH formula and see that,

\begin{equation}
I_{\Delta t} u(\cdot,0)-u(\cdot,\Delta t)=\frac{1}{2}\Delta t^2\big([\mathcal{L}_3,\mathcal{L}_2]+[\mathcal{L}_3,\mathcal{L}_1]+[\mathcal{L}_2,\mathcal{L}_1]\big)u(\cdot,0)+\mathcal{O}(\Delta t^3).
\end{equation} Zeros in $\mathcal{O}(1)$ and $\mathcal{O}(\Delta t)$ term show that the splitting scheme is locally consistent, which can be equivalently {achieved} by series expansion in terms of $\Delta t$. Moreover, we find that computing the $k$-th order modified equation associated with Eq.\eqref{eqn:particleSDE} in BEA is equivalent to computing the terms of BCH formula up to order $(\Delta t)^k$ 
in the Eq.\eqref{eqn:operatorflow}. 
{We can see that the solution generated by Eq.\eqref{scheme} follows a perturbed Hamiltonian system (with divergence-free velocity and additive noise) at any order ${k}$, by considering the $({k}+1)$-nested Lie bracket consisting of $\{-f\partial_{x_1},\ g\partial_{x_2},\ \partial^{2}_{x_1x_1}+\partial^{2}_{x_2x_2}\}$. Moreover, we can easily derive that they {generate divergence-free} fields.}

In \cite{WangXinZhang:18}, we proved that for the SDE \eqref{eqn:particleSDE} with a time-dependent and separable Hamiltonian $H(x_1,x_2,t)=F({x_2},t)+G({x_1},t)$,  the numerical solution obtained by using the symplectic-preserving scheme \eqref{scheme} follows an asymptotic Hamiltonian $H^{\Delta t}(x_1,x_2,t)$, which is a first-order approximation to 
$H(x_1,x_2,t)$. Equivalently, the velocity field in the first-order modified backward Kolmogorov equation is divergence-free and the invariant measure on the torus (defined by $\mathbb{R}^d/\mathbb{Z}^d$, when period is $1$) remains uniform, which is also known as the Haar measure. However, the numerical solution obtained using the Euler-Maruyama scheme for the SDE \eqref{eqn:particleSDE} does not have these properties. 

Moreover, given any explicit splitting scheme for deterministic systems, by adding additive noise we shall have a similar form of flow propagation. And we shall see in later proof that, such operator formulation is very effective in analyzing the order of convergence and volume-preserving property.


\section{Convergence analysis}\label{sec:ConvergenceAnalysis}
\noindent 
We shall prove the convergence rate of our symplectic stochastic integrators in computing effective diffusivity based on a probabilistic approach, which allows us to get rid of the exponential growth factor in our error estimate. {As stated at the beginning of Section \ref{sec:DerivationSchemes}, we will first limit our analysis to 2D separable Hamiltonian velocity fields. We will show in Section \ref{sec:HighDimensionalCases} that all the derivations can be generalized to high-dimensional cases.}
\subsection{Convergence to an invariant measure}\label{sec:InvariantMeasure}
\noindent 
The numerical method to compute effective diffusivity of a passive tracer model is closely related to study the limit of a sequence generated by the stochastic integrators. Therefore, we can apply the 
results from ergodic theory to study the convergence of the solution. 
The following result is fundamental for the proof of our convergence analysis. 
\begin{proposition}\label{pro:doobs}
On the torus space $\tilde{\textbf{Y}}=\mathbb{R}^{2}/\mathbb{Z}^{{2}}$, let $I_{\Delta t}^*$ denote the  transform of the density function during $\Delta t$ using the numerical scheme \eqref{scheme}. 	Let $I_{\Delta t}$ denote the adjoint operator (i.e., the flow operator) of $I_{\Delta t}^*$ in the space of $\mathcal{B}(\tilde{\textbf{Y}})$, which is the set of bounded measurable functions on $\tilde{\textbf{Y}}$. 
	Then, $I_{\Delta t}$ is a compact operator from $\mathcal{B}(\tilde{\textbf{Y}})$ to itself. And there exists one and only one invariant probability measure on $(\tilde{Y},\Sigma)$, denoted as $\pi$, satisfying,
	
	\begin{equation}
	\sup_{x\in\tilde{Y}}\Big|(I^n_{\Delta t}\phi) (\textbf{x})-\int \phi(\textbf{x}')\pi(d\textbf{x}')\Big|\leq C ||\phi||_{L_\infty}e^{-\rho n},\quad \forall \phi\in \mathcal{B}(\tilde{\textbf{Y}}),
	\end{equation}
	where $\rho>0$, $C>0$ are independent {of} $\phi(\cdot)$.
	\begin{proof}
		We shall verify that the transition kernel associated with the numerical scheme \eqref{scheme} 
		satisfies the assumptions required by the Theorem 3.3.1 (see the page 199 in \cite{BensoussanLionsPapa:2011}).
		First in the $\mathbb{R}^2$ space, the integration process associated with the numerical scheme 
		can be expressed as a Markov process with the transition kernel,
		
		\begin{align}\label{eqn:rnkernel}
		K_{\Delta t}&\big((x_1^{n-1},x_2^{n-1}),(x_1^{n},x_2^{n})\big)=\nonumber\\&\frac{1}{2\pi \sigma^2\Delta t}\exp\Bigg(
		-\frac{\Big(x_1^{n}-x_1^{n-1}+f(x_2^{n-1})\Delta t\Big)^2+\Big(x_2^{n}-x_2^{n-1}-g\big(x_1^{n-1}-f(x_2^{n-1})\Delta t\big)\Delta t\Big)^2}{2\sigma^2\Delta t}\Bigg),
		\end{align}	
		where$(x_1^{n},x_2^{n})$ is the solution obtained by applying the scheme \eqref{scheme} on $(x_1^{n-1},x_2^{n-1})$ with time step $\Delta t$. 
		
		{ Since $f$ and $g$ are periodic functions, we can project the solution of SDE \eqref{eqn:particleSDE} on the torus space $\tilde{\textbf{Y}}=\mathbb R^2/\mathbb Z^2$ pathwisely. We denote the solution on the torus as $\tilde{\textbf{X}}$ and its numerical approximation as $\tilde{\textbf{X}^n}$. Given any periodic function $f$, we know $f(\textbf{X})=f|_{\tilde{\textbf{Y}}}(\tilde{\textbf{X}})$. Later on, for simplicity reasons, we do not distinguish $f$ and $f|_{\tilde{\textbf{Y}}}$. Moreover, we do not distinguish $\textbf{X}$ and $\tilde{\textbf{X}}$ when 
		we apply a periodic function on it.  Eq.\eqref{eqn:rnkernel} can be directly extended to the torus space $\tilde{\textbf{Y}}$ as } 
		
		\begin{align}\label{eqn:tnkernel}
		\tilde{K}_{\Delta t}&\big((x_1^{n-1},x_2^{n-1}),(x_1^{n},x_2^{n})\big)=\sum_{i,j\in\mathbb{Z}}\frac{1}{2\pi \sigma^2\Delta t}\cdot\nonumber\\&\exp\Bigg(-\frac{\Big(x_1^{n}+i-x_1^{n-1}+f(x_2^{n-1})\Delta t\Big)^2+\Big(x_2^{n}+j-x_2^{n-1}-g\big(x_1^{n-1}-f(x_2^{n-1})\Delta t\big)\Delta t\Big)^2}{2\sigma^2\Delta t}\Bigg).
		\end{align}
		One can see that if $0 < \Delta t\ll1$, then $\tilde{K}$ is smooth and is essentially bounded above zero, 
		i.e., $essn~\tilde{K}>0, ~\forall \big((x_1^{n-1},x_2^{n-1}),(x_1^{n},x_2^{n})\big)\in \tilde{Y}\times\tilde{Y}$. 
		Thus, the operator $I_{\Delta t}$ is compact since it is an integral operator with a smooth kernel. Then applying the Theorem 3.3.1 in \cite{BensoussanLionsPapa:2011}, we prove the assertion of the Proposition \ref{pro:doobs}. 
	\end{proof}
\end{proposition}
Now, we state a corollary that is a simple conclusion of exponential decay property proved in Proposition \ref{pro:doobs},
which will be useful in the proof of main results of this paper.
\begin{corollary}\label{col:ergodicity}
	Given that the assumptions in Proposition \ref{pro:doobs} are satisfied and $\phi\in\mathcal{B}(\tilde{Y})$, we have for all initial $\textbf{X}^0\in{\mathbb{R}^2}$
	
	\begin{equation}
	\lim_{n\to\infty}\frac{1}{n}\sum_{i=1}^n \mathbb{E} \phi(\textbf{X}^i)=\int_{\tilde{Y}} 
	\phi(\textbf{x})\pi(d\textbf{x}).
	\end{equation}
\end{corollary}

Before we close this subsection, we present a convergence result for the inverse of operator sequences, which can also be viewed as a modification of Theorem 1.16 in Section IV of \cite{kato2013perturbation}.
\begin{proposition}\label{prop:inverseoperator}

Let $\mathcal{X},\mathcal{Y}$ denote two Banach spaces. Assume $T_{n }$, $T$ are bounded linear operators from  $\mathcal{X}$ to $\mathcal{Y}$, satisfying $\lim_{n \to \infty}||T_{n }-T||_{\mathcal{B}(\mathcal{X},\mathcal{Y})}=0$, and $T^{-1}\in\mathcal{B}(\mathcal{Y},\mathcal{X})$. Given $f\in \mathcal{Y}$, if $T_{n}^{-1}f$, $n=1,2,...$ uniquely exist, then we have a convergence estimate as follows,

	\begin{align}\label{eqn:est_inverseoperator}
	\lim_{n \to \infty}\big|\big|(T_{n}^{-1}-T^{-1})f\big|\big|=0.
	\end{align}
	\begin{proof}
		After some simple calculations, we get	
		
		\begin{align}
		T_{n}^{-1}-T^{-1}&=T^{-1}(T-T_{n})T_{n}^{-1} \nonumber\\
		&=T^{-1}(T-T_{n})T^{-1}+T^{-1}(T-T_{n})(T_{n}^{-1}-T^{-1}).\label{eqn:inverse_indentity}
		\end{align}
		Now applying $T_{n}^{-1}-T^{-1}$ on $f$, we get
		
		\begin{align}
		||(T_{n}^{-1}-T^{-1})f||\leq & ||T^{-1}||^2\cdot||T-T_{n}||\cdot||f||\nonumber\\
		&+ ||T^{-1}||\cdot||T-T_{n}||\cdot ||(T_{n}^{-1}-T^{-1})f||
		\end{align}
		Since $\lim_{n \to \infty}||T_{n}-T||=0$, we assume for $n\geq N_0$, $||T_{n}-T||\cdot ||T^{-1}||<\frac{1}{2}$, then,
		
		\begin{align}
		||(T_{n}^{-1}-T^{-1})f||\leq 2||T^{-1}||^2\cdot||T-T_{n}||\cdot||f||,\quad \forall n\geq N_0,
		\label{eqn:est_inverseoperator0}
		\end{align}
		Eq.\eqref{eqn:est_inverseoperator} follows if we take the limit as $n \to \infty$ on both sides of
		\eqref{eqn:est_inverseoperator0}. 
	\end{proof}
\end{proposition}
\subsection{A discrete-type cell problem}\label{sec:ProbabilisticProof}
\noindent 
In the Eulerian framework, the periodic solution of the cell problem \eqref{CellProblem_EffectiveDiffusivity} and the corresponding formula for the effective diffusivity \eqref{Def_EffectiveDiffusivity_Euler} play a key role in studying the behaviors of the chaotic and stochastic flows. In the Lagrangian framework,
we shall define a discrete analogue of the cell problem that enables us to compute the effective diffusivity. We revisit the scheme Eq.\eqref{scheme},  

\begin{equation}\label{scheme_analysis}
\begin{cases}
x_{1}^{n}=x_{1}^{n-1}-f(x_{2}^{n-1})\Delta t +\sigma N^{n-1}_{x_1}\\
x_{2}^{n}=x_{2}^{n-1}+g\big(x_{1}^{n-1}-f(x_{2}^{n-1})\Delta t\big)\Delta t +\sigma N^{n-1}_{x_2},
\end{cases}
\end{equation}
where $N^{n-1}_{x_1}$, $N^{n-1}_{x_2}\sim\sqrt{\Delta t}\mathcal{N}(0,1)$ are i.i.d. normal random variables.

We will show that the solutions $x_{1}^{n}$ and $x_{2}^{n}$ obtained by the scheme \eqref{scheme_analysis} have bounded expectations if the initial values are bounded. Taking expectation of the first equation of Eq.\eqref{scheme_analysis} on both sides, we obtain 

\begin{align}
\mathbb{E}x_{1}^{n}=\mathbb{E}x_{1}^{n-1}-\Delta t \mathbb{E}f(x_{2}^{n-1})=\mathbb{E}x_{1}^{0}-\Delta t\sum_{k=0}^{n-1}\mathbb{E}f(x_{2}^{k}).
\label{bounded_Ep}
\end{align}
{ 
As a symplectic scheme in 2D, \eqref{scheme_analysis} admits the uniform measure as its invariant measure. Then applying Proposition \ref{pro:doobs} and using the fact that $f$ is a periodic function with zero mean, we know that,}

\begin{equation}\label{eqn:decay_expactation}
\sup_{(x_{1}^{0},x_{2}^{0})\in{\mathbb{R}^2}}\big|\mathbb{E}f(x_{2}^{k})\big|\leq e^{-\rho k}||f||_\infty.
\end{equation}
By applying triangle inequalities in Eq.\eqref{bounded_Ep} and using the result in  Eq.\eqref{eqn:decay_expactation}, we arrive at,

\begin{equation}
|\mathbb{E}x_{1}^{n}|\leq|\mathbb{E}x_{1}^{0}|+C_1||f||_\infty, \label{bounded_Ep2}
\end{equation}
where $C_1$ does not depend on $n$. Using the same approach, we know that $\mathbb{E}x_{2}^{n}$ is also bounded. Now, we are in the position to define the discrete-type cell problem. { Recalling that $\textbf{X}^n=(x_{1}^{n},x_{2}^{n})^T$ denotes the solution of discrete scheme at $t_n=n\Delta t$, we first define} 

\begin{equation}\label{eqn:cell_prob_f}
\hat{f}(\textbf{x})=-\Delta t\sum_{n=0}^{\infty}\mathbb{E}[f(\textbf{X}^n)|\textbf{X}^0=\textbf{x}],\quad \textbf{x}\in{\mathbb{R}^2},
\end{equation}
where the summability is guaranteed by Eq.\eqref{eqn:decay_expactation}. {{$f(\textbf{X}^n)$ is equivalent to $f(x_{2}^{n})$ in our case. This is due to that the velocity fields are given by separable Hamiltonian functions, so $f(\textbf{X}^n)=f(x_{1}^{n},x_{2}^{n})$ is independent of $x_{1}^{n}$.} At the same time, we should notice that $\hat{f}(\textbf{x})$ relies on the second component of $\textbf{x}$, as the initial condition is $\textbf{X}^0=\textbf{x}$}. Then, we shall show that $\hat{f}(\textbf{x})$ satisfies the following properties. 
\begin{lemma}\label{lem:existenceofcellpro}
According to our assumption on the Hamiltonian, which is separable and periodic along each dimension, we know that $f$ is a periodic function with zero mean on $\tilde{Y}$, i.e., $\int_{\tilde{Y}}f=0$. Therefore, $\hat{f}$ defined in \eqref{eqn:cell_prob_f} is the unique solution in $\mathcal{B}_0(\tilde{Y})$ such that,

\begin{equation}
\hat{f}(\textbf{X}^0)+\Delta tf(\textbf{X}^0)= \mathbb{E}[\hat{f}(\textbf{X}^{1})|\textbf{X}^0].
\end{equation} 
Moreover, {$\hat{f}$} is smooth.
\begin{proof}

{ Starting from Eq.\eqref{eqn:cell_prob_f} and by the periodicity of $f$, we know that $\hat{f}$ is a periodic function. }	{Then, by using basic properties of
	conditional expectation, we can get that}

		\begin{align}
		\hat{f}(\textbf{X}^0)+\Delta tf(\textbf{X}^0)=&\Delta t\mathbb{E}[\sum_{m=0}^{\infty}-f(\textbf{X}^m)|\textbf{X}^0]+\Delta t f(\textbf{X}^0) = -\Delta t\mathbb{E}[\sum_{m=1}^{\infty}f(\textbf{X}^m)|\textbf{X}^0]\nonumber\\
		=&-\Delta t\mathbb{E}\big[\mathbb{E}[\sum_{m=1}^{\infty}f(\textbf{X}^m)|\textbf{X}^1]|\textbf{X}^0\big] = \mathbb{E}[\hat{f}(\textbf{X}^{1})|\textbf{X}^0].
		\label{DiscreteCellProblem1}
		\end{align}	
		Recall the definition of the operator \eqref{eqn:operatorflow}, Eq.\eqref{DiscreteCellProblem1}
		implies that  
		\begin{equation}\label{cell_eqn}
		(I_{\Delta t}-I_d)\hat{f}= I_{\Delta t}\hat{f}-\hat{f}=\Delta tf,
		\end{equation}
where $I_d$ is the identity operator. Moreover, since $f$ is smooth and the mapping of the operator $I_{\Delta t}$ on bounded functions will generate smooth functions, so $\hat{f}$ is smooth. 
		
{According to Proposition \ref{pro:doobs}, the invariant (measure) of $I_{\Delta t}^*$ is unique and it is the uniform measure. In other words, the null space of the operator $I_{\Delta t}^*-I_d$ consists of constant functions. Then following the assumption that $f$ is mean zero on $\tilde{\textbf{Y}}$, we know $f$ is in $\mathcal{N}(I_{\Delta t}^*-I_d)^\perp$.   By the Fredholm alternative with the fact that $I_{\Delta t}$ is a compact operator, we arrive at the conclusion that the solution $\hat{f}$ to Eq.\eqref{cell_eqn} is unique in $\mathcal{B}(\tilde{\textbf{Y}})$ up to a constant and 
		it smoothly depends on $f$. }
	\end{proof}
\end{lemma}
{Noticing that the passive tracer model \eqref{eqn:particleSDE} is autonomous}, we obtain 
\begin{equation}\label{key_eqn}
\mathbb{E}[\hat{f}(\textbf{X}^{n+1})|\textbf{X}^n]-\hat{f}(\textbf{X}^n)=\Delta tf(\textbf{X}^n),\quad a.s.\quad\forall n\in \mathbb{N}.
\end{equation}
\begin{remark}
{For the second component of the solution $\textbf{X}^n$, i.e., $x_{2}^{n}$, we can define the discrete cell problem in the same manner. Notice the numerical schemes for $x_{1}^{n}$ and $x_{2}^{n}$ have the same structures}. As such, we define
	\begin{equation}
	\hat{g}(\textbf{x})=\Delta t\sum_{n=0}^{\infty}\mathbb{E}[g(\textbf{X}^{',n})|\textbf{X}^0=\textbf{x}],\quad \textbf{x}\in{\mathbb{R}^2},
	\end{equation}
where {$\textbf{X}^{',n}=\textbf{X}^n-\Delta t \left(f(\textbf{X}^n),0\right)^T$. {Under the assumption that the drift terms $f$ and $g$ in Eq.\eqref{eqn:particleSDE} are smooth, we know the leading order term of $g(\textbf{X}^{',n})$ is $g(\textbf{X}^n)$. Then, we can carry out the analysis for $\hat{g}(\textbf{x})$ in the same manner as that for $\hat{f}(\textbf{x})$.}}

\end{remark} 

The Proposition \ref{pro:doobs} and the Lemma \ref{lem:existenceofcellpro} are very general results. 
In the remaining part of this paper, we only need the result that $\hat{f}$ is unique in an H\"{o}lder space $\mathbb{C}^{p,\alpha}_0(\tilde{\textbf{Y}})\subsetneq\mathcal{B}(\tilde{\textbf{Y}})$. To be precise, given a smooth drift function
$f$, $\hat{f}$ shall be in $\mathbb{C}^{p,\alpha}_0(\tilde{Y})$, where $p\geq 6, 0<\alpha<1$ and the subscript 
index $0$ indicates that it is a subspace with zero-mean functions. To prove that  $I_{\Delta t}$ is a compact operator from $\mathbb{C}^{p,\alpha}_0(\tilde{Y})$ to itself is quite standard. We can apply the Arzel\`{a}-Ascoli theorem to verify the relative compactness of the operator $I_{\Delta t}$ by studying its mapped results on a bounded set. Both equicontinuity and  point-wise boundedness come as the result that $I_{\Delta t}$ is an integral operator with a smooth kernel.    
However, we do not want to complicate the presentation by pursuing this avenue.

\subsection{Convergence estimate of the discrete-type cell problem}\label{sec:ConvergenceDiscreteProblem}
\noindent
After defining the discrete-type cell problem (e.g., Eq.\eqref{cell_eqn}) and proving the existence and uniqueness of the solution $\hat{f}$, we shall prove that $\hat{f}$ converges to the solution of a continuous cell problem in certain subspace, e.g., $\mathbb{C}^{6,\alpha}_0(\tilde{\textbf{Y}})$. We remark that in the remaining part of this paper, we shall choose the space $\mathbb{C}^{6,\alpha}_0(\tilde{\textbf{Y}})$ to carry out our analysis. However there is no requirement that we have to choose this space. In fact, any space that has certain regularity (belongs to the domain of the operator $\mathcal{L}$) will work. To start with, we define the following continuous cell problem 
\begin{equation}\label{continuous_cellproblem}
\mathcal{L}\chi_1=f, 
\end{equation}
where the operator $\mathcal{L}$ is defined in Eq.\eqref{HamiltonianFlowOperator}. 
Given $f$ is a  smooth function defined on $\tilde{\textbf{Y}}$ with zero mean, the Eq.\eqref{continuous_cellproblem} admits a unique solution $\chi_1$ in $\mathbb{C}_0^{6,\alpha}(\tilde{\textbf{\textbf{Y}}})$. This is a standard result of elliptic PDEs in H\"{o}lder space (see, e.g., the Theorem 6.5.3 in \cite{krylov1996lectures}). Moreover, $\mathcal{L}$ is a bijection between two Banach spaces $\mathbb{C}_0^{6,\alpha}(\tilde{\textbf{Y}})$ and $\mathbb{C}_0^{4,\alpha}(\tilde{\textbf{Y}})$, and its inverse is bounded. The following theorem
states that under certain conditions the solution of the discrete-type cell problem converges
to the solution of the continuous one.
\begin{theorem} \label{est:cellproblem}
Assume $f$ is a smooth function defined on $\tilde{\textbf{Y}}$ with zero mean. 
Let $\hat{f}$ and $\chi_1$ be the solutions to the discrete-type cell problem \eqref{cell_eqn} 
and continuous cell problem \eqref{continuous_cellproblem}, respectively.  
When $\Delta t\to 0$, the solution $\hat{f}$ converges to the solution $\chi_1$ in $\mathbb{C}^{p,\alpha}_0$, at the rate of $\mathcal{O}(\Delta t)$,  where $p\geq 6$ and $0<\alpha<1$.
\end{theorem}
\begin{proof}
Integrating Eq.\eqref{continuous_cellproblem} along time gives, 

\begin{align}\label{eqn:asymcell}
\exp (\Delta t \mathcal{L})\chi_1-\chi_1=f\Delta t +\mathcal{O}((\Delta t)^2):=\Delta t\bar{f},
\end{align}
where $\bar{f}=f+O(\Delta t)$. Combining Eqns.\eqref{cell_eqn} and \eqref{eqn:asymcell}, we obtain
\begin{equation}\label{eqn:derivationform}
\exp (\Delta t \mathcal{L})\chi_1 - I_{\Delta t} \hat{f} 
-(\chi_1 - \hat{f}  ) = \Delta t(\bar{f}-f)
\end{equation}
Eq.\eqref{eqn:derivationform} shows the connection between $\chi_1$ and $\hat{f}$. 
After some simple calculations, we get
\begin{equation}\label{eqn:balancedform}
\mathcal{L}(\chi_1-\hat{f})=(\mathcal{L}-\tilde{L}_1)(\chi_1-\hat{f})+\tilde{L}_2\hat{f}+(\bar{f}-f),
\end{equation}
where
\begin{equation}
\tilde{L}_1:=\frac{\exp (\Delta t \mathcal{L})-I_d}{\Delta t},
\quad\text{and} \quad 
\tilde{L}_2:=\frac{I_{\Delta t}-\exp (\Delta t \mathcal{L})}{\Delta t}.
\end{equation}
One can easily verify that in the space of bounded linear operators from $\mathbb{C}_0^{6,\alpha}(\tilde{\textbf{Y}})$ to $\mathbb{C}_0^{4,\alpha}(\tilde{\textbf{Y}})$, there is a strong convergence in the operator norm $||\cdot||$, 

\begin{equation}\label{conv:operatorstrong}
||\tilde{L}_1-\mathcal{L}||=\mathcal{O}(\Delta t)\quad \text{as } \Delta t\to 0.
\end{equation}
For the operator $\tilde{L}_2$, by using the BCH formula \eqref{eqn:BCHformula} we can obtain, 

\begin{align}\label{convergenceofhatf2}
\tilde{L}_2\to& \frac{\exp\Big(\frac{\Delta t^2}{2}\big([L_3,L_2]+[L_2,L_1]+[L_3,L_1]\big)+\mathcal{O}((\Delta t)^3)\Big)-I_{d}}{\Delta t}\cdot \exp (\Delta t \mathcal{L}) \nonumber\\
\to& \frac{\Delta t}{2}\big([L_3,L_2]+[L_2,L_1]+[L_3,L_1]\big) +\mathcal{O}((\Delta t)^2).
\end{align}
Denoting $\tilde{L}_3:=\tilde{L}_1+\tilde{L}_2\equiv \frac{I_{\Delta t}-I_d}{\Delta t}$, we have  $\tilde{L}_3\to \mathcal{L}$ in $\mathcal{B}\big(\mathbb{C}_0^{6,\alpha}(\tilde{\textbf{Y}}),\mathbb{C}_0^{4,\alpha}(\tilde{\textbf{Y}})\big)$. 
Finally, applying the Proposition \ref{prop:inverseoperator}, we get,
\begin{equation}\label{convergenceofhatf}
\lim_{\Delta t \to 0}\hat{f}=\lim_{\Delta t \to 0}\tilde{L}_3^{-1}f=\mathcal{L}^{-1}f=\chi_1.
\end{equation}
In addition, combining the results of the Eqns.\eqref{eqn:asymcell}, \eqref{conv:operatorstrong}, \eqref{convergenceofhatf2} and \eqref{convergenceofhatf} for the right hand side of Eq.\eqref{eqn:balancedform}, we know that when $\Delta t$ { is }small enough (does not depend on the total computational time $T$, but may depend on the estimate of $f$, $g$ and $\sigma$), the following convergence estimate holds
\begin{equation}\label{ConvergenceResult_ftochi}
||\chi_1-\hat{f}||=\mathcal{O}(\Delta t).
\end{equation}
Thus, the assertion in Theorem \ref{est:cellproblem} is proved.   
\end{proof}
\subsection{Convergence estimate for the effective diffusivity}\label{sec:EstimateEffectiveDiffusivity}
\noindent
We shall show the main estimates in this section. We first prove that the second-order moment
of the solution obtained by using our numerical scheme has an (at most) linear growth rate. Secondly, we provide the convergence rate of our method in computing the effective diffusivity.
\begin{theorem}\label{thm:boundness}
Let {$\textbf{X}^n=(x_1^n,x_2^n)^T$} denote the solution of the passive tracer model \eqref{eqn:particleSDE}
obtained by using our numerical scheme with time-step $\Delta t$.  If the Hamiltonian $H(x_1,x_2)$ is separable, periodic and smooth enough (in order to guarantee the existence and uniqueness of the solution to the SDE \eqref{eqn:particleSDE}), then we can prove that the second-order moment of the solution ${\textbf{X}^n}$ (a discrete Markov process) is at most linear growth, i.e.,
\begin{equation}\label{conj}
\max_n\big\{\mathbb{E}\frac{||{\textbf{X}^n}||^2}{n}\big\}\ \text{is bounded.}
\end{equation}
\begin{proof} 
		We first estimate the second-order moment of the first component of ${\textbf{X}^n}=(x_1^n,x_2^n)^T$, since the other one can be estimated in the same manner. Simple calculations show that 
		
		\begin{align}
		\mathbb{E}[(x_1^n)^2|(x_1^{n-1},x_2^{n-1})]&=\mathbb{E}\big(x_1^{n-1}-f(x_2^{n-1})\Delta t+\sigma N^{n-1}_{x_1} \big)^2\nonumber\\
		&=\mathbb{E}(x_1^{n-1})^2+\Delta t \big(\sigma^2-2\mathbb{E}[x_1^{n-1}f(x_2^{n-1})]\big)+
		(\Delta t)^2\mathbb{E}(f(x_2^{n-1}))^2.
		\end{align}
		We should point out that the term $\mathbb{E}[x_1^{n-1}f(x_2^{n-1})]$ corresponds to the convection enhanced level of the diffusivity. Our goal is to prove that the term $\mathbb{E}[x_1^{n-1}f(x_2^{n-1})]$ is bounded over $n$, though it may depend on 
		$f$, $g$ and $\sigma$. {To be noted that, here we are calculating the expectation of $(x_1^n)^2$, which is not defined in the torus space. But in the following derivation we will show that it can be decomposed into sums of periodic functions acting on $\textbf{X}^n=(x_1^n,x_2^n)^T$. Hence after the decomposition (see Eq.\eqref{eqn:var}) we can still apply the previous analysis on torus space.}  
		
		We now directly compute the contribution of the term $\mathbb{E}[x_1^{n-1}f(x_2^{n-1})]$ to the effective diffusivity with
		the help of Eq.\eqref{key_eqn},  
		
		\begin{align}\label{ComputeEpfq}
		\Delta t\sum_{i=0}^{n-1}\mathbb{E}[x_1^if(x_2^i)]=\sum_{i=0}^{n-1}\mathbb{E}\big[x_1^i \big(\mathbb{E}[\hat{f}(\textbf{X}^{i+1})|\textbf{X}^i]-\hat{f}(\textbf{X}^i)\big)\big].
		\end{align}
		 {Throughout the proof, 
		we shall use the fact that if $\mathbf{X}$, $\mathbf{Y}$ are random processes and $\mathbf{Y}$ is measurable under a filtration $\mathcal{F}$, then with appropriate integrability assumption, we have 
		
		\begin{equation}
		\mathbb{E}[\mathbf{X}\mathbf{Y}]=\mathbb{E}\Big[\mathbb{E}[\mathbf{X}\mathbf{Y}|\mathcal{F}]\Big]
		=\mathbb{E}\Big[\mathbb{E}[\mathbf{X}|\mathcal{F}]\mathbf{Y}\Big].
		\end{equation}} Let $\mathcal{F}_i$ denote the filtration generated by the solution process until $\textbf{X}^i$. Notice that $x_1^i\in\mathcal{F}_i$, for the Eq.\eqref{ComputeEpfq}, we have
	
		\begin{align}
		RHS&=\sum_{i=0}^{n-1}\mathbb{E}\big[x_1^i \big(\hat{f}(\textbf{X}^{i+1})-\hat{f}(\textbf{X}^i)\big)\big]
		\nonumber\\
		&=\sum_{i=1}^{n}\mathbb{E}\big[\hat{f}(\textbf{X}^i)(x_1^{i-1}-x_1^i)\big]-\hat{f}(\textbf{X}^0)x_1^0+\mathbb{E}[\hat{f}(\textbf{X}^n)x_1^n]\nonumber\\
		&= \sum_{i=1}^{n}\mathbb{E}\big[\hat{f}(\textbf{X}^i)\big(f(x_1^{i-1})\Delta t-\sigma N^{i-1}_{x_1}\big)\big]-\hat{f}(\textbf{X}^0)x_1^0+\mathbb{E}[\hat{f}(\textbf{X}^n)x_1^n].
		\end{align}
		Hence,
		
		\begin{align}\label{eqn:var}
		\frac{1}{n}\mathbb{E}\big[(x_1^n)^2|(x_1^{0},x_2^{0})\big]=&\frac{1}{n}(x_1^0)^2+\Delta t\sigma^2-2\Delta t\frac{1}{n}\sum_{i=0}^{n-1}\mathbb{E}[x_1^if(x_2^i)]+(\Delta t)^2\frac{1}{n}\sum_{i=0}^{n-1}\mathbb{E}f^2(x_2^i)\nonumber\\
		=&\frac{1}{n}(x_1^0)^2+\Delta t\sigma^2+(\Delta t)^2\frac{1}{n}\sum_{i=0}^{n-1}\mathbb{E}f^2(x_2^i)
		-\frac{2}{n}\sum_{i=1}^{n}\mathbb{E}\big[\hat{f}(\textbf{X}^i)\big(f(x_2^{i-1})\Delta t-\sigma N^{i-1}_{x_1}\big)\big]\nonumber\\
		&-\frac{2}{n}\big(\hat{f}(\textbf{X}^0)x_1^0-\mathbb{E}[\hat{f}(\textbf{X}^n)x_1^n]\big).
		\end{align}
		Recall the fact that $\textbf{X}^n=(x_1^n,x_2^n)$ converges to the uniform measure in distribution. So given any continuous periodic function $f^*$, the Corollary \ref{col:ergodicity} implies
		
		\begin{equation}
		\lim_{n\to\infty} \mathbb{E}f^*(\textbf{X}^n)=\int_{\tilde{Y}}f^*(\textbf{x})d\textbf{x}.
		\end{equation}
		Furthermore, we have the estimate 
		
		\begin{equation}\label{est:limsup}
		\limsup_{n\to\infty}  {\mathbb{E}}\frac{1}{n}\sum_{i=0}^nf^*(\textbf{X}^i)<\infty.
		\end{equation}
		Applying the Cauchy-Schwarz inequality for the term $\frac{2}{n}\sum_{i=1}^{n}\mathbb{E}\big[\hat{f}(\textbf{X}^i)\big(f(x_2^{i-1})\Delta t-\sigma N^{i-1}_{x_1}\big)\big]$
		in Eq.\eqref{eqn:var} and replacing $f^*$ by $f^2$ and $\hat{f}^2$ in Eq.\eqref{est:limsup}, we can prove that $\frac{1}{n}\mathbb{E}\big[(x_1^n)^2|(x_1^{0},x_2^{0})\big]$ is bounded.  Using the same trick, we know that $\frac{1}{n}\mathbb{E}\big[(x_2^n)^2|(x_1^{0},x_2^{0})\big]$ is also bounded. Thus, the assertion in Eq.\eqref{conj} is proved. 
	\end{proof}
\end{theorem}		
In our numerical scheme \eqref{scheme}, we first fix the time-step $\Delta t$ and use it to compute the effective diffusivity until the result converges to a constant, which may depend on $\Delta t$. Next, we shall prove that the limit of the constant converges to the exact effective diffusivity of the original passive tracer model as $\Delta t$ approaches zero. Namely, we shall prove that our numerical scheme is robust in computing the effective diffusivity.
\begin{theorem}\label{thm:convergence}
	Let $x_1^n$, $n=0,1,....$ be the numerical solution of the first component of the scheme \eqref{scheme} and  $\Delta t$ denote the time-step. We have the convergence estimate of the effective diffusivity as 
	
	\begin{align}\label{est:order}
	\lim_{n\to\infty}\frac{\mathbb{E}(x_1^n)^2}{n\Delta t}= \sigma^2-2\int_{\mathbb{T}^2} \chi_1f+\mathcal{O}(\Delta t),
	\end{align}
	where the constant in $\mathcal{O}(\Delta t)$ does not depends on the computational time $T$.
	\begin{proof}
		We divide both sides of the Eq.\eqref{eqn:var} by $\Delta t$ and obtain 
		
		\begin{align}\label{eqn:mainterm}
		\frac{1}{n\Delta t}\mathbb{E}[(x_1^n)^2|(x_1^{0},x_2^{0})]=&\frac{1}{n\Delta t}(x_1^0)^2+\sigma^2+ \frac{\Delta t}{n}\sum_{i=0}^{n-1}\mathbb{E}f^2(x_2^i)\nonumber\\
		&-\frac{2}{n\Delta t}\sum_{i=1}^{n}\mathbb{E}\big[\hat{f}(\textbf{X}^i)\big(f(x_2^{i-1})\Delta t-\sigma N^{i-1}_{x_1}\big)\big]\nonumber\\
		&-\frac{2}{n\Delta t}\big(\hat{f}(\textbf{X}^0)x_1^0-\mathbb{E}[\hat{f}(\textbf{X}^n)x_1^n]\big)
		\end{align}
		First, we notice that for a fixed $\Delta t$, the terms $\frac{1}{n\Delta t}(x_1^0)^2$ and 
		$\frac{2}{n\Delta t}\hat{f}(\textbf{X}^0)x_1^0$ converge to zero as $n\to\infty$, where we have used the fact $\hat{f}(\textbf{X}^0)$ is bounded. Then, for a fixed $\Delta t$, we have 
		
		\begin{equation}
		\lim_{n\to\infty}\frac{2}{n\Delta t}\big|\mathbb{E}[\hat{f}(\textbf{X}^n)x_1^n]\big|\leq \lim_{n\to\infty}\frac{2}{\sqrt{n}\Delta t}||\hat{f}||_{\infty}\mathbb{E}|\frac{x_1^n}{\sqrt{n}}|\leq \lim_{n\to\infty}\frac{1}{\sqrt{n}\Delta t}||\hat{f}||_{\infty}\mathbb{E}[\frac{(x_1^n)^2}{n}+1]=0,
		\end{equation}
		where the term $\mathbb{E}[\frac{(x_1^n)^2}{n}]$ is bounded due to the Theorem \ref{thm:boundness} and $||\hat{f}||_\infty\to||\chi_1||_\infty<\infty $ due to the Theorem \ref{est:cellproblem}.
		Therefore, we only need to focus on the estimate of terms in the second line of Eq.\eqref{eqn:mainterm}, which correspond to the convection-enhanced diffusion effect. 
		Notice that $\hat{f}\in\mathbb{C}^{6,\alpha}$, we compute the Ito-Taylor series approximation of $\hat{f}(\textbf{X}^i)$,
		
		\begin{align}\label{eqn:taylor}
		\hat{f}(\textbf{X}^i)=&\hat{f}(\textbf{X}^{i-1})+\hat{f}_{x_1}(\textbf{X}^{i-1})\big(-f(x_2^{i-1})\Delta t+\sigma N^{i-1}_{x_1}\big)+\hat{f}_{x_2}(\textbf{X}^{i-1})\big(g(x_1^{i-1})\Delta t+\sigma N^{i-1}_{x_2}\big)\nonumber\\
		&+\frac{1}{2}\big(\hat{f}_{x_1x_1}(X^{i-1})+\hat{f}_{x_2x_2}(\textbf{X}^{i-1})\big)\sigma^2\Delta t+\mathcal{O}(\Delta t^2).
		\end{align}
		Since $\hat{f}\to\chi_1$ in $\mathbb{C}_0^{6,\alpha}$, the truncated term $\mathcal{O}(\Delta t^2)$ in Eq.\eqref{eqn:taylor} is uniformly bounded when $\Delta t$ is small enough. Substituting the Taylor 
		expansion of $\hat{f}(\textbf{X}^i)$ into the target term of our estimate, we get 
		
		\begin{align}
		\mathbb{E}[\hat{f}(\textbf{X}^i)&(f(x_2^{i-1})\Delta t-\sigma N^{i-1}_{x_1})]=\mathbb{E}\Big[\Big(f(x_2^{i-1})\Delta t-\sigma N^{i-1}_{x_1}\Big)\cdot\nonumber\\
		&\Big(\hat{f}(\textbf{X}^{i-1})+\hat{f}_{x_1}(\textbf{X}^{i-1})\big(-f(x_2^{i-1})\Delta t+\sigma N^{i-1}_{x_1}\big)\nonumber\\
		&+\hat{f}_{x_2}(\textbf{X}^{i-1})\big(g(x_1^{i-1})\Delta t+\sigma N^{i-1}_{x_2}\big)+\frac{1}{2}\big(\hat{f}_{x_1x_1}(\textbf{X}^{i-1})+\hat{f}_{x_2x_2}
		(\textbf{X}^{i-1})\big)\sigma^2\Delta t+\mathcal{O}(\Delta t^2)\Big)\Big].
		\end{align}
		Combining the terms with the same order of $\Delta t$, we obtain 
		
		\begin{align}
		\mathbb{E}\big[\hat{f}(\textbf{X}^i)\big(f(x_2^{i-1})\Delta t-\sigma N^{i-1}_{x_1}\big)\big]=\Delta t\mathbb{E}[\hat{f}(\textbf{X}^{i-1})f(x_2^{i-1})-\sigma^2\hat{f}_{x_1}(\textbf{X}^{i-1})]+\mathcal{O}(\Delta t^2),
		\end{align}
where we have used the facts that: (1) $\textbf{X}^{i-1}$ is independent {of} $N^{i-1}_{x_1}$ and $N^{i-1}_{x_2}$ so the expectations of the corresponding terms vanish; (2) $N^{i-1}_{x_1}$ and $N^{i-1}_{x_2}$ are independent so $\mathbb{E}N^{i-1}_{x_1}N^{i-1}_{x_2}=0$; and (3) $\mathbb{E}(N^{i-1}_{x_1})^2=\Delta t$. Finally, by using the Corollary \ref{col:ergodicity} and noticing the invariant measure is the uniform measure, we obtain from Eq.\eqref{eqn:mainterm} that

		\begin{align}
		\lim_{n\to\infty} \frac{1}{n\Delta t}\mathbb{E}[(x_1^n)^2|(x_1^{0},x_2^{0})]=\sigma^2-2\int (\hat{f}f-\sigma^2\hat{f}_{x_1})+\mathcal{O}(\Delta t).
		\end{align}
		Thus, our statement in the Eq.\eqref{est:order} is proved 
		using the facts that $\hat{f}$ converges to $\chi_1$ (see Theorem \ref{est:cellproblem})
		and $\int \hat{f}_{x_1}=0$.
	\end{proof}
\end{theorem}
\begin{remark}
If we divide two on both sides of the Eq.\eqref{est:order}, we can find that our result 
recovers the definition of the effective diffusivity $D^{E}_{11}$ defined in the Eq.\eqref{Def_EffectiveDiffusivity_Euler}. This reveals the connection of the definition of the effective diffusivity using the Eulerian framework and Lagrangian framework.	
\end{remark}
	
\subsection{Generalizations to high-dimensional cases}\label{sec:HighDimensionalCases}  
\noindent
To show the essential idea of our probabilistic approach, we have carried out our convergence analysis based on a two-dimensional model problem \eqref{eqn:particleSDE}. In fact, the extension of our approach to higher-dimensional problems is straightforward. Now we consider a high-dimensional problem as follow, 
\begin{equation}
\label{eqn:multiDSDE}
d\textbf{X}(t)=\textbf{v}(\textbf{X}(t))dt+\Sigma d\textbf{W}(t),
\end{equation}
where $\textbf{X}=(x_1,x_2,\cdots,x_d)^{T}\in \mathbb R^{d}$ is the position of a particle, 
$\textbf{v}=(v_1,v_2,\cdots,v_d)^{T}\in \mathbb R^{d}$ is the Eulerian velocity field at position $X$, $\Sigma$ is a $d\times d$ constant non-singular matrix, and $\textbf{W}(t)$ is a $d$-dimensional Brownian motion vector. In particular, we assume the $v_i$ does not depend on $x_i$, $i=1,...,d$.  Thus, the incompressible condition for $\textbf{v}(\textbf{X})$ (i.e. $\nabla_{\textbf{X}}\cdot \textbf{v}(\textbf{X})=0$) is easily guaranteed. 

For a deterministic and divergence-free dynamical system, Feng et. al. proposed a volume-preserving method \cite{KangShang1995volume}, which splits a $d$-dimensional problem into $d-1$ subproblems with each of them being a two-dimensional problem and thus being volume-preserving. We shall modify Feng's method (first-order case) by including the randomness as the last subproblem to take into account the additive noise, i.e., 
\begin{equation}\label{scheme:multiDimension}
	{\begin{cases}
			x_1^{*}=x_1^{n-1}+\Delta t v_1(x_2^{n-1},x_3^{n-1},x_4^{n-1},\cdots,x_{d-1}^{n-1}, x_d^{n-1}),\\
			x_2^{*}=x_2^{n-1}+\Delta t v_2(x_1^{*},x_3^{n-1},x_4^{n-1},\cdots,x_{d-1}^{n-1}, x_d^{n-1}),\\
			x_3^{*}=x_3^{n-1}+\Delta t v_3(x_1^{*},x_2^{*},x_4^{n-1},\cdots,x_{d-1}^{n-1}, x_d^{n-1}),\\
			\cdots,\\ 
			x_d^{*}=x_d^{n-1}+\Delta t v_d(x_1^{*},x_2^{*},x_3^{*},x_4^{*},\cdots, x_{d-1}^{*}),\\
			\textbf{X}^{n}=\textbf{X}^{*}+\Sigma(\textbf{W}^{n}-\textbf{W}^{n-1}),
	\end{cases}}
\end{equation}
where $\textbf{X}^{*}=(x_1^{*},x_2^{*},\cdots,x_d^{*})^{T}$, $\textbf{W}^{n}-\textbf{W}^{n-1}$ is a $d$-dimensional independent random vector with each component of the form $\sqrt{\Delta t}\xi_{i}$, $\xi_{i}\sim\mathcal{N}(0,1)$, and $\textbf{X}^{n}=(x_1^{n},x_2^{n},\cdot\cdot\cdot,x_d^{n})^T$  is the numerical approximation to the exact solution $\textbf{X}(t_n)$ to the SDE \eqref{eqn:multiDSDE}  at time $t_n=n\Delta t$. 

The techniques of the convergence analysis for two-dimensional problem can be applied to 
high-dimensional problems without much difficulty. For the high-dimensional problem \eqref{eqn:multiDSDE}, 
the smoothness and  strict positivity of the transition kernel in the discrete process can be
guaranteed if one assumes that the covariance matrix $\Sigma$ is non-singular and the scheme 
\eqref{scheme:multiDimension} is explicit. According to our assumption for the velocity field, 
the scheme \eqref{scheme:multiDimension} is volume-preserving. Thus, the solution to the first-order modified equation is divergence-free and the invariant measure on the torus (defined by $\mathbb{R}^d/\mathbb{Z}^d$, when period is $1$) remains uniform. Finally, the convergence of the cell problem can be studied by using the BCH formula \eqref{eqn:BCHformula} with $d+1$ 
PDE operators. Recall that in the Eq.\eqref{eqn:operatorflow} we have three PDE operators when we 
study the two-dimensional problem. Therefore, our numerical methods are robust in computing effective diffusivity for high-dimensional problems, which will be demonstrated through the three-dimensional chaotic flow problems in the Section \ref{sec:Numerical_example}. 
\section{Numerical Examples}\label{sec:Numerical_example}
\noindent
The aim of this section is two-fold. First, we shall design challenging numerical examples to verify the convergence analysis proposed in this paper, especially the Theorem \ref{thm:convergence}. Secondly, we shall investigate the diffusion enhancement for several chaotic  velocity fields. Without loss of generality, we compute the quantity $\frac{\mathbb{E}[x_1(T)^2]}{2T}$, which is used to approximate $D^{E}_{11}$ in the effective diffusivity matrix \eqref{Def_EffectiveDiffusivity_Euler}. 

\subsection{Verification of the convergence rate}\label{sec:VerifyConvergenceRate}
\noindent
We first consider a passive tracer model, where the velocity field is given by a chaotic cellular flow with oscillating vortices. Specifically, the flow is generated by a Hamiltonian defined as 

\begin{align}
H(x_1,x_2)=\frac{1}{2\pi}\exp(\sin(2\pi x_1))-\frac{1}{4\pi}\exp(\cos(4\pi x_2+1)).
\end{align}
The motion of a particle moving in this chaotic cellular flow is described by the SDE,
\begin{equation}\label{eqn:CellularFlowOscillatingVortices}
\begin{cases}
dx_1=\sin(4\pi x_2 +1)\exp(\cos(4\pi x_2+1))dt+\sigma dW_1,\\
dx_2= \cos(2\pi x_1)\exp(\sin(2\pi x_1)) dt+\sigma dW_2,
\end{cases}
\end{equation}
where $\sigma=\sqrt{2\times 0.01}$, $W_{i}$ are independent Brownian motions, and the initial data $(x_1^0,x_2^0)$ follows uniform distributions in $[-0.5,0.5]^2$. 

In our numerical experiments, we use Monte Carlo samples to discretize the Brownian motions $W_1$ and $W_2$. The sample number is denoted by $N_{mc}$. We choose $\Delta t_{ref}=0.001$ and $N_{mc}=640,000$ to solve the SDE \eqref{eqn:CellularFlowOscillatingVortices} and compute the reference solution, i.e., the ``exact'' effective diffusivity, where the final computational time is $T=12000$ so that the calculated effective diffusivity converges to a constant. It takes about 20 hours to compute the reference solution on a 64-core server (Gridpoint System at HKU). The reference solution for the effective diffusivity is $D^{E}_{11}=0.12629$.  

\begin{figure*}[htbp]
	\centering
	\begin{subfigure}{0.5\linewidth}
		\includegraphics[width=\linewidth]{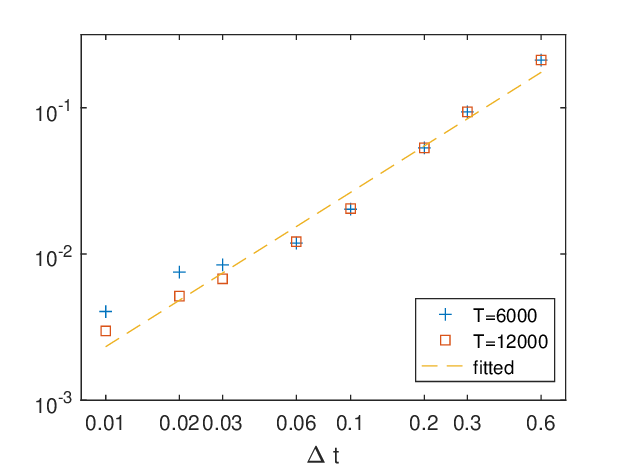}
		\caption{2D chaotic cellular flow, fitted slope $\approx 1.04$}
		\label{fig:meg1f1new}
	\end{subfigure}%
	~
	\begin{subfigure}{0.5\linewidth}
		\includegraphics[width=\linewidth]{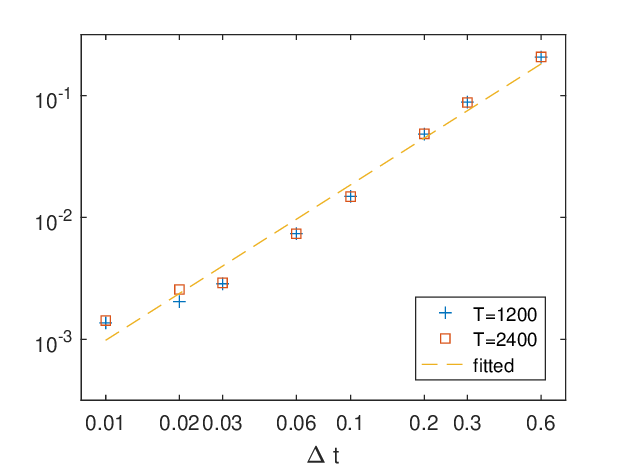}
		
		\caption{3D Kolmogorov-type flow, fitted slope $\approx 1.27$}
		\label{fig:meg2f1new}
	\end{subfigure}
	\caption{Error of $D^{E}_{11}$ in different computational times and flows with different time-steps.}
	
\end{figure*}

In Fig.\ref{fig:meg1f1new}, we plot the convergence results of the 
effective diffusivity using our method (i.e., $\frac{\mathbb{E}[x_1(T)^2]}{2T}$) with respective to different time-step $\Delta t$ at $T=6000$ and $T=12000$. The computational time of our method depends on 
$N_{mc}$, $\Delta t$, and $T$. In this example, it takes less than two hours to get the one associated with the $N_{mc}=640,000$, $\Delta t=0.01$, and $T=12000$. In addition, we show a fitted straight line with the slope $1.04$, i.e., the convergence rate is about $(\Delta t)^{1.04}$. Meanwhile, by comparing two sets of data in the Fig.\ref{fig:meg1f1new}, corresponding to the numerical effective diffusivity obtained at different computational times, we can see that error does not grow with respect to time, which justifies the statement in Theorem  \ref{thm:convergence}. 

To further study the accuracy and robustness of our numerical method in solving high-dimensional problems, we consider a 3D Kolmogorov-type flow. Let $(x_1,x_2,x_3)^T\in R^3$ denote the position of a particle in the 3D Cartesian coordinate system. The motion of a particle moving in the 3D Kolmogorov-type flow is described by the following SDE,
\begin{equation}\label{eqn:3Dflow}
\begin{cases}
dx_1= \cos(4\pi x_3+1)\exp(\sin(4\pi x_3+1))dt+\sigma dW_1,\\
dx_2= \cos(6\pi x_1+2)\exp(\sin(6\pi x_1+2))dt+\sigma dW_2,\\
dx_3= \cos(2\pi x_2+3)\exp(\sin(2\pi x_2+3))dt+\sigma dW_3,\\
\end{cases}
\end{equation} 
where $W_{i}$ are independent Brownian motions. This is inspired by the so-called Kolmogorov flow \cite{KflowGalloway:1992} (see Eq.\eqref{eqn:Kflow}). The Kolmogorov flow is obtained from the Arnold-Beltrami-Childress (ABC) flow with $A=B=C=1$ and with cosines taken out. Behaviors of the classic Kolmogorov flow will be discussed later.

In our numerical experiments, 
we choose $\Delta t_{ref}=0.001$ and $N_{mc}=6,400,000$ to solve the SDE \eqref{eqn:3Dflow} and compute the reference solution, i.e., the ``exact'' effective diffusivity.  After some numerical tests, we find that the passive tracer model will enter a mixing stage if the computational time is set to be $T=2400$. It takes about 56 hours to compute the reference solution on the server and the reference solution for the effective diffusivity is 
$D^{E}_{11}=0.13106$. 

In Fig. \ref{fig:meg2f1new}, we plot the convergence results of the 
effective diffusivity using our method with respect to different time-step $\Delta t$. In addition, we show a fitted straight line with the slope $1.27$, i.e., the convergence rate is about $(\Delta t)^{1.27}$. This numerical result also agrees with our error analysis.

\subsection{Investigation of the convection-enhanced diffusion phenomenon}\label{subsec:ExperimentOnPhenomenon}
\begin{figure}[tbph]
	\centering
	\includegraphics[width=0.6\linewidth]{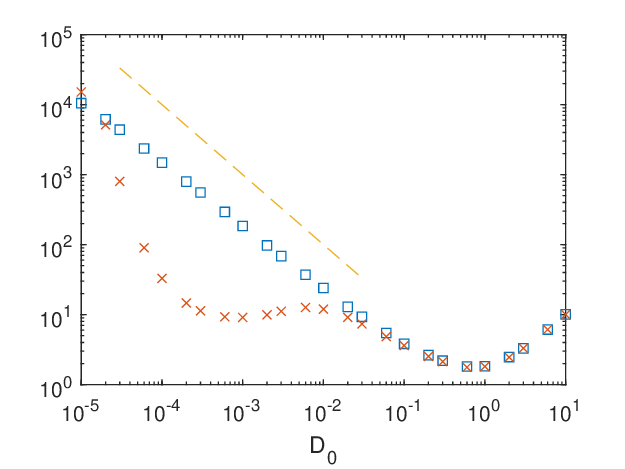}
	\caption{Convection-enhanced diffusion with maximal enhancement in ABC flow: $\Box$ for the symplectic scheme, $\times$ for the Euler scheme, $--$ for reference line $y=\frac{1}{D_0}$.}
	\label{fig:eg7f1}
\end{figure}

\noindent 
We first consider the classical ABC flow with our symplectic stochastic integrators.
The ABC flow is a three-dimensional incompressible velocity field which is an exact solution to the Euler's equation. It is notable as a simple example of a fluid flow that can have chaotic trajectories. 
The particle is transported by the velocity field $v=(A\sin(x_3)+C\cos(x_2),B\sin(x_1)+A\cos(x_3),C\sin(x_2)+B\cos(x_1))$ and perturbed by an additive noise. The associated passive tracer model reads 
\begin{equation}\label{eqn:ABCflow}
\begin{cases}
dx_1= (A\sin(x_3)+C\cos(x_2))dt+\sigma dW_1,\\
dx_2= (B\sin(x_1)+A\cos(x_3))dt+\sigma dW_2,\\
dx_3= (C\sin(x_2)+B\cos(x_1))dt+\sigma dW_3,\\
\end{cases}
\end{equation}
where $W_{i}$ are independent Brownian motions. In Fig.\ref{fig:eg7f1}, we show the relation between $D^E_{11}$ and $D_0$. Recall that the parameter $D_0=\sigma^2/2$. By setting $A=B=C=1$, we recover the same phenomenon as the Fig.2 in \cite{Biferale:95}, for $D_0\in[10^{-3},10^{-1}]$ and can extend to $D_0\in[10^{-5},10^{-4}]$
; see Fig.\ref{fig:eg7f1}. { As a comparison to our stochastic structure-preserving scheme, we directly apply the Euler-Maruyama scheme (also called the Euler scheme) to solve the SDE \eqref{eqn:ABCflow}. We can see that the Euler scheme failed {to recover it} when $D_0$ is small. The evidence for the failure of the Euler scheme when $D_0$ is small can be also found in \cite{WangXinZhang:18}.  } 
The Fig.\ref{fig:eg7f1} shows that the $D^E_{11}$ of the ABC flow obtained by our symplectic method corresponds to upper-bound of Eq.\eqref{eqn:maximaldiffusion}, i.e. the maximal enhancement, $D^E_{11}\sim \mathcal{O}(1/D_0)$. This maximal enhancement phenomenon may be attributed to the ballistic  orbits of the ABC flow, which was discussed in \cite{mcmillen2016ballistic,xin2016periodic}. 

\begin{figure}[t!]
	\centering
	\begin{subfigure}[t]{0.5\textwidth}
		\centering
		\includegraphics[width=\linewidth]{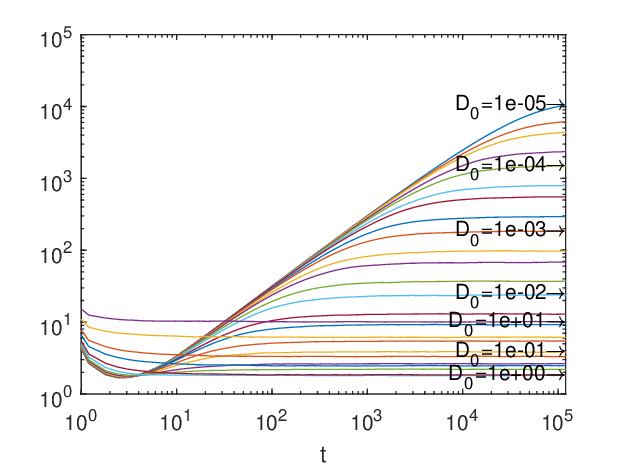}
		\caption{$\frac{\mathbb{E}[x_1(t)^2]}{2t}$ of different $D_0$ in the symplectic scheme}
		\label{fig:eg7f2}
	\end{subfigure}%
	~ 
	\begin{subfigure}[t]{0.5\textwidth}
		\centering
		\includegraphics[width=\linewidth]{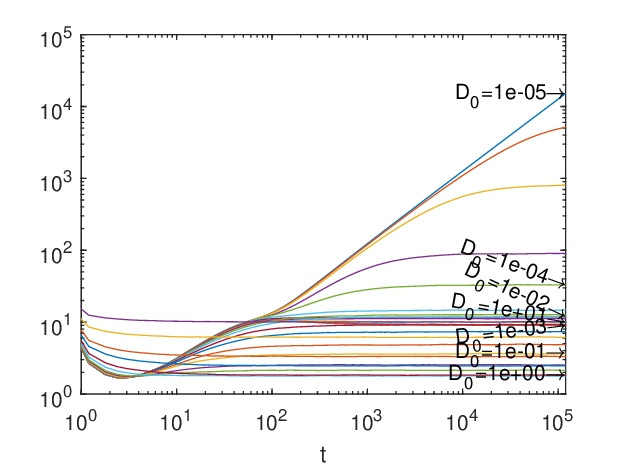}
		\caption{$\frac{\mathbb{E}[x_1(t)^2]}{2t}$ of different $D_0$ in the Euler scheme}
		\label{fig:eg7f3}
	\end{subfigure}
	\caption{Calculated $D^E_{11}$ in the ABC flow along time via two different schemes}
\end{figure}

From Fig.\ref{fig:eg7f2} we can see that diffusion time, i.e., the time when $\frac{\mathbb{E}[x_1(t)^2]}{2t}$ approaches a constant, increases as $\mathcal{O}(1/D_0)$ when $D_0\to 0$ in the symplectic scheme. 
{Interested readers are referred to \cite{feng2019dissipation} to find that the upper bound of diffusion time can be a bit smaller than $\mathcal{O}(1/D_0)$ given the strong mixing property of the flows.
Due to the gap between chaotic and strongly mixing flows, to the best of our knowledge,  the diffusion time (as $D_0$ tends to $0$) for chaotic flows has yet to be rigorously proved. Fig.\ref{fig:eg7f2} shows the diffusion time of ABC flow may reach the upper bound in the a priori estimate for general flows. However, the Euler scheme gives a different result in  Fig.\ref{fig:eg7f3}. It attains a diffusion time which is much faster than $O(1/D_0)$. This may { be} due to the numerical dissipation of the Euler scheme.} The statement that the Euler scheme generates wrong results can also be found in the Fig.\ref{fig:eg7f1}.   

\begin{figure}[htbp]
	\centering
	\includegraphics[width=0.6\linewidth]{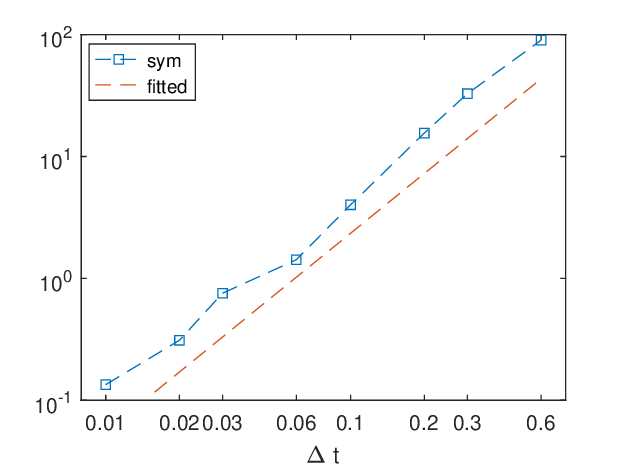}
	\caption{Error of $D_{11}^E$ in the ABC flow, the dashed line with $\Box$ is for the symplectic scheme, 
		and the slope of the fitted is $\approx 1.42$.} 
	\label{fig:eg7f4}
\end{figure}

We point out that the error estimate in  Theorem \ref{thm:convergence} is just an upper bound. Fig.\ref{fig:eg7f4} shows that when $D_0$ is $10^{-3}$, the convergence rate is about $\mathcal{O}(\Delta t^{1.42})$. It is very expensive to study the passive tracer model for the ABC flow since the diffusing time is extremely long. In our numerical test for the  Fig.\ref{fig:eg7f4}, we choose $N_{mc}=120,000$, $\Delta t=0.001$, and $T=12,000$. In this setting, the error of the Monte Carlo simulation cannot be avoided, so there is a small oscillation around the fitted slope. 

\begin{figure*}[bthp]
	\centering
	\includegraphics[width=0.6\linewidth]{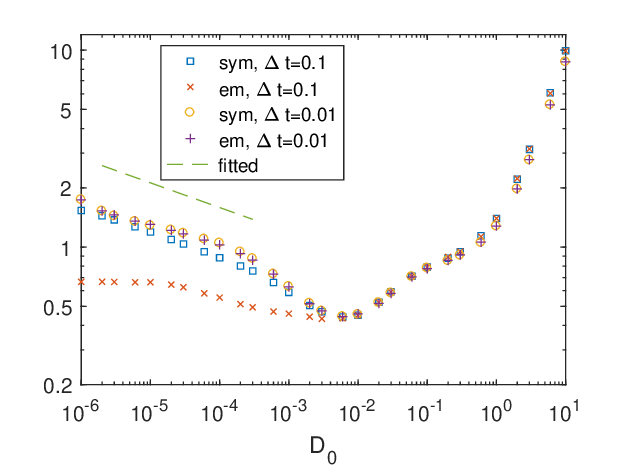}
	\caption{Convection-enhanced diffusion with sub-maximal enhancement in Kolmogorov flow. ``sym'' means the results for symplectic scheme and ``em'' means the results for Euler scheme. $--$ means the fitted line for small $D_0$ with slope $\approx -0.13$.}
	\label{fig:eg8f1}
\end{figure*}

Finally, we investigate the  convection-enhanced diffusion phenomenon for another chaotic flow, i.e., the Kolmogorov flow. The associated passive tracer model reads,
\begin{equation}\label{eqn:Kflow}
\begin{cases}
dx_1= \sin(x_3)dt+\sigma dW_1,\\
dx_2= \sin(x_1)dt+\sigma dW_2,\\
dx_3= \sin(x_2)dt+\sigma dW_3,\\
\end{cases}
\end{equation} 
where $W_{i}$ are independent Brownian motions. In Fig.\ref{fig:eg8f1}, we show the relation between $D^E_{11}$ and $D_0$, where $D_0=\sigma^2/2$. For each $D_0$, we use $N_{mc}=120,000$ particles to solve the SDE \eqref{eqn:Kflow} via the symplectic method and the Euler method with $\Delta t=0.1$ . The final computational time is $T=12,000$ so that the particles are fully mixed for $D_0\geq10^{-6}$. 

Under such setting, we find that the dependency of $D^E_{11}$ on $D_0$ is quite different from the chaotic and stochastic flows that we have studied in \cite{WangXinZhang:18} and from the foregoing ABC flow (maximal enhancement). The fitted slope within $D_0\in[10^{-6},10^{-5}]$ is $-0.13$, which indicates that $D^E_{11}\sim \mathcal{O}(1/D_0^{0.13})$. {The slope is significantly greater than $-1$ and this can be called sub-maximal enhancement. The existence of sub-maximal enhancement may be explained by the fact that the Kolmogorov flow is more chaotic than the ABC flow \cite{KflowGalloway:1992}. }The chaotic trajectories in Kolmogorov flow enhance diffusion much less than channel like structures such as the ballistic orbits of ABC flows \cite{mcmillen2016ballistic,xin2016periodic}.  More studies on the diffusion enhancement phenomenon of the ABC flow and the Kolmogorov flow, especially the time-dependent cases 
will be reported in our future work. 

We also compare the performance of the symplectic scheme and Euler scheme in computing the 
effective diffusivity for the Kolmogorov flow. Specifically, we implement the symplectic scheme and Euler scheme with 
time step  $\Delta t=0.1$ and $\Delta t=0.01$, respectively. 
In Fig.\ref{fig:eg8f1}, we find that (1) the symplectic scheme with $\Delta t=0.1$ and $\Delta t=0.01$ will give similar 
results in computing the effective diffusivity; (2) the symplectic scheme and the Euler scheme with $\Delta t=0.01$ 
will give almost the same convergent results in computing the effective diffusivity, which provides evidence that 
our statement on the Kolmogorov flow (i.e., the sub-maximal enhancement phenomenon) is correct; (3) the Euler scheme with $\Delta t=0.1$ gives wrong results but the symplectic scheme with $\Delta t=0.1$ gives acceptable results, which 
provides evidence that the symplectic scheme is very robust in computing the effective diffusivity. In this example, the symplectic scheme approximately achieves a $10\times$ speedup over the Euler scheme.

\begin{figure*}[htbp]
	\centering
	\begin{subfigure}[t]{0.5\textwidth}
		\centering
		\includegraphics[width=\linewidth]{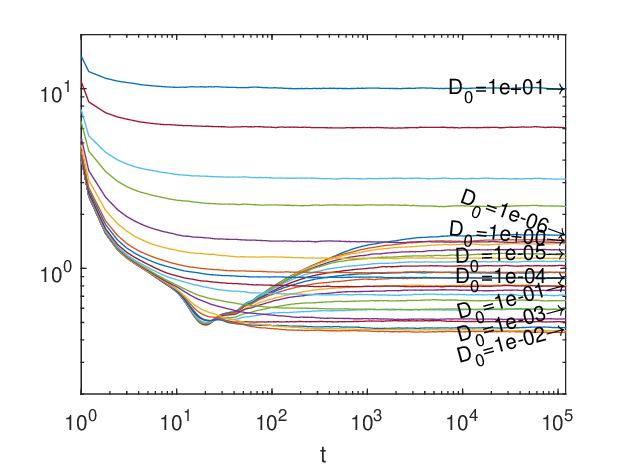}
		\caption{$\frac{\mathbb{E}[x_1(t)^2]}{2t}$ of different $D_0$ using the symplectic scheme}
		\label{fig:eg8f2}
	\end{subfigure}%
	~ 
	\begin{subfigure}[t]{0.5\textwidth}
		\centering
		\includegraphics[width=\linewidth]{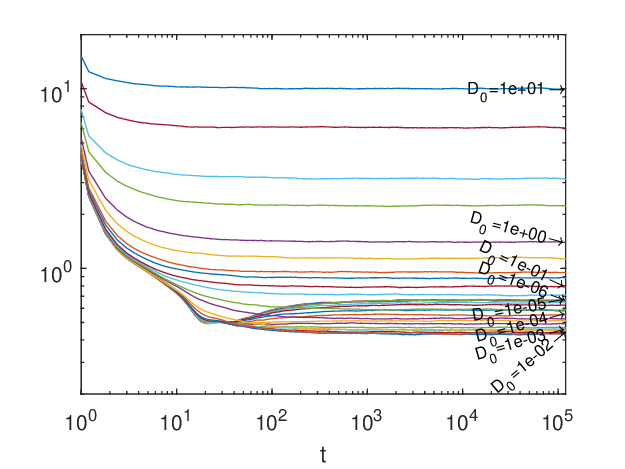}
		\caption{$\frac{\mathbb{E}[x_1(t)^2]}{2t}$ of different $D_0$ in Euler scheme}
		\label{fig:eg8f3}
	\end{subfigure}
	\caption{Calculated $D^E_{11}$ in the Kolmogorov flow via two different schemes.}
\end{figure*}

\begin{figure}[hbtp]
	\centering
	\includegraphics[width=0.6\linewidth]{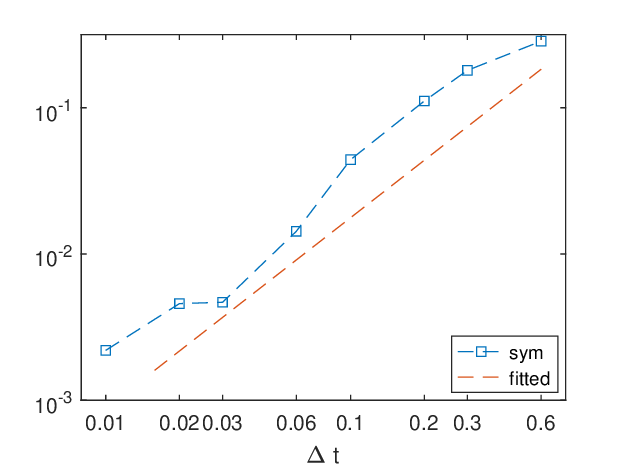}
	\caption{Error of $D_{11}^E$ in the Kolmogorov flow. The slope of the fitted line is $\approx 1.30$.}
	\label{fig:eg8f4}
\end{figure}

Fig.\ref{fig:eg8f2} and Fig.\ref{fig:eg8f3} show different behaviors of the numerical effective diffusivity $\frac{\mathbb{E}[x_1(t)^2]}{2t}$ obtained using the symplectic scheme and the Euler scheme with respect to computational time. Specifically, Fig.\ref{fig:eg8f2} shows $T=12000$ is quite enough for $D_0\geq10^{-6}$. And in Fig.\ref{fig:eg8f3}, it seems that in Euler scheme, the diffusion time is much smaller. {
Similar to our investigation in ABC flows, this may { be} due to the excess numerical dissipation generated by the Euler scheme.}	
In Fig.\ref{fig:eg8f4}, we also study the convergence rate of the symplectic scheme in computing the
effective diffusivity for the Kolmogorov flow \eqref{eqn:Kflow}. We find that the convergence rate is {  $\mathcal{O}(\Delta t^{1.3})$} in this example. 

\section{Conclusions}\label{sec:conclusion}
\noindent  
In this paper, we analyzed the robustness of a numerical scheme to compute the effective diffusivity of passive tracer models, especially for the three-dimensional ABC flow and the Kolmogorov flow.  
The scheme is based on the Lagrangian formulation of the passive tracer model, i.e., solving SDEs. We split the SDE problem into a deterministic sub-problem and a stochastic one, where the former is discretized using a symplectic-preserving scheme while the later is solved using the {Euler} scheme.  We provide a completely new error analysis for our numerical scheme that is based on a probabilistic approach, which gives a sharp and uniform in time
error estimate for the numerical solution of the effective diffusivity. Finally, we present numerical results to demonstrate the accuracy  of the proposed method for several typical chaotic flow problems of physical interests, including the Arnold-Beltrami-Childress (ABC) flow and the Kolmogorov flow. We observed the maximal enhancement phenomenon in the ABC flows and the sub-maximal enhancement phenomenon in the Kolmogorov flow, respectively.  

There are two directions we plan to explore in our future work. First, we shall extend the probabilistic approach to provide sharp convergence analysis in computing effective diffusivity for time-dependent chaotic flows, such as time-dependent ABC flows. In addition, we shall investigate the convection-enhanced diffusion phenomenon for general spatial-temporal stochastic flows \cite{Yaulandim:1998,Majda:99} and develop convergence analysis for the corresponding numerical methods.

\section*{Acknowledgments}
\noindent
The research of Z. Wang is partially supported by the Hong Kong PhD Fellowship Scheme. The research of J. Xin is partially supported by NSF grants DMS-1211179, DMS-1522383 and IIS-1632935. The research of Z. Zhang is supported by Hong Kong RGC grants (Projects 27300616, 17300817, and 17300318), National Natural Science Foundation of China (Project 11601457), Seed Funding Programme for Basic Research (HKU), and an RAE Improvement Fund from the Faculty of Science (HKU). The computations were performed using the HKU ITS research computing facilities that are supported in part by the Hong Kong UGC Special Equipment Grant (SEG HKU09).

\bibliographystyle{plain}
\bibliography{ZWpaper}

\begin{thebibliography}{10}

\bibitem{abdulle2015long}
A.~Abdulle, G.~Vilmart, and K.~Zygalakis.
\newblock Long time accuracy of {L}ie--{T}rotter splitting methods for
  {L}angevin dynamics.
\newblock {\em SIAM Journal on Numerical Analysis}, 53(1):1--16, 2015.

\bibitem{JanHesthaven2017structure}
B.~Afkham and J.~Hesthaven.
\newblock Structure preserving model reduction of parametric hamiltonian
  systems.
\newblock {\em SIAM Journal on Scientific Computing}, 39(6):A2616--A2644, 2017.

\bibitem{BenOwhadi2003}
G.~Ben~Arous and H.~Owhadi.
\newblock Multiscale homogenization with bounded ratios and anomalous slow
  diffusion.
\newblock {\em Communications on Pure and Applied Mathematics}, 56(1):80--113,
  2003.

\bibitem{BensoussanLionsPapa:2011}
A.~Bensoussan, J.~L. Lions, and G.~Papanicolaou.
\newblock {\em Asymptotic analysis for periodic structures}, volume 374.
\newblock American Mathematical Soc., 2011.

\bibitem{Biferale:95}
L.~Biferale, A.~Crisanti, M.~Vergassola, and A.~Vulpiani.
\newblock Eddy diffusivities in scalar transport.
\newblock {\em Phys. Fluids}, 7:2725--2734, 1995.

\bibitem{bou2010long}
N.~Bou-Rabee and H.~Owhadi.
\newblock Long-run accuracy of variational integrators in the stochastic
  context.
\newblock {\em SIAM Journal on Numerical Analysis}, 48(1):278--297, 2010.

\bibitem{Carmona1997homogenization}
R.~Carmona and L.~Xu.
\newblock Homogenization for time-dependent two-dimensional incompressible
  {G}aussian flows.
\newblock {\em The Annals of Applied Probability}, 7(1):265--279, 1997.

\bibitem{Fannjiang:94}
A.~Fannjiang and G.~Papanicolaou.
\newblock Convection-enhanced diffusion for periodic flows.
\newblock {\em SIAM J Appl. Math.}, 54:333--408, 1994.

\bibitem{Fannjiang:97}
A.~Fannjiang and G.~Papanicolaou.
\newblock Convection-enhanced diffusion for random flows.
\newblock {\em J. Stat. Phys.}, 88:1033--1076, 1997.

\bibitem{KangShang1995volume}
K.~Feng and Z.~Shang.
\newblock Volume-preserving algorithms for source-free dynamical systems.
\newblock {\em Numerische Mathematik}, 71(4):451--463, 1995.

\bibitem{feng2019dissipation}
Y.~Feng and G.~Iyer.
\newblock Dissipation enhancement by mixing.
\newblock {\em Nonlinearity}, 32(5):1810, 2019.

\bibitem{KflowGalloway:1992}
D.~Galloway and M.~Proctor.
\newblock Numerical calculations of fast dynamos in smooth velocity fields with
  realistic diffusion.
\newblock {\em Nature}, 356(6371):691, 1992.

\bibitem{Garnier:97}
J.~Garnier.
\newblock Homogenization in a periodic and time-dependent potential.
\newblock {\em SIAM Journal on Applied Mathematics}, 57(1):95--111, 1997.

\bibitem{BCHformula1974baker}
R.~Gilmore.
\newblock Baker-{C}ampbell-{H}ausdorff formulas.
\newblock {\em Journal of Mathematical Physics}, 15(12):2090--2092, 1974.

\bibitem{ErnstLubich:06}
E.~Hairer, C.~Lubich, and G~Wanner.
\newblock {\em Geometric numerical integration: structure-preserving algorithms
  for ordinary differential equations}.
\newblock Springer Science and Business Media, 2006.

\bibitem{hong2006multi}
J.~Hong, H.~Liu, and G.~Sun.
\newblock The multi-symplecticity of partitioned runge-kutta methods for
  hamiltonian pdes.
\newblock {\em Mathematics of computation}, 75(253):167--181, 2006.

\bibitem{Oleinik:94}
V.~V. Jikov, S.~Kozlov, and O.~A. Oleinik.
\newblock {\em Homogenization of {D}ifferential {O}perators and {I}ntegral
  {F}unctionals}.
\newblock Springer, Berlin, 1994.

\bibitem{kato2013perturbation}
Tosio Kato.
\newblock {\em Perturbation theory for linear operators}, volume 132.
\newblock Springer Science \& Business Media, 2013.

\bibitem{krylov1996lectures}
N.~V. Krylov.
\newblock {\em Lectures on elliptic and parabolic equations in H{\"o}lder
  spaces}.
\newblock Graduate studies in mathematics.

\bibitem{Yaulandim:1998}
C.~Landim, S.~Olla, and H.~T. Yau.
\newblock Convection--diffusion equation with space--time ergodic random flow.
\newblock {\em Probability theory and related fields}, 112(2):203--220, 1998.

\bibitem{JackXin:11}
Y.~Liu, J.~Xin, and Y.~Yu.
\newblock Asymptotics for turbulent flame speeds of the viscous {G}-equation
  enhanced by cellular and shear flows.
\newblock {\em Arch. Rational Mech. Anal.}, 202:461--492, 2011.

\bibitem{JackXinLyu:2017}
J.~Lyu, J.~Xin, and Y.~Yu.
\newblock Computing residual diffusivity by adaptive basis learning via
  spectral method.
\newblock {\em Numerical Mathematics: Theory, Methods and Applications},
  10(2):351--372, 2017.

\bibitem{Majda:99}
A.~J. Majda and P.~R. Kramer.
\newblock Simplified models for turbulent diffusion: theory, numerical
  modelling, and physical phenomena.
\newblock {\em Phys. Rep.}, 314:237--574, 1999.

\bibitem{mclaughlin1997effect}
R.~McLaughlin and J.~Zhu.
\newblock The effect of finite front thickness on the enhanced speed of
  propagation.
\newblock {\em Combustion science and technology}, 129(1-6):89--112, 1997.

\bibitem{mcmillen2016ballistic}
T.~McMillen, J.~Xin, Y.~F. Yu, and A.~Zlatos.
\newblock Ballistic orbits and front speed enhancement for abc flows.
\newblock {\em SIAM Journal on Applied Dynamical Systems}, 15(3):1753--1782,
  2016.

\bibitem{mezic1996maximal}
I.~Mezi{\'c}, J.~F. Brady, and S.~Wiggins.
\newblock Maximal effective diffusivity for time-periodic incompressible fluid
  flows.
\newblock {\em SIAM Journal on Applied Mathematics}, 56(1):40--56, 1996.

\bibitem{Milstein:02}
G.~Milstein, Y.~Repin, and M.~Tretyakov.
\newblock Symplectic integration of {H}amiltonian systems with additive noise.
\newblock {\em SIAM J. Numer. Anal}, 39:2066--2088, 2002.

\bibitem{Oksendal:13}
B.~Oksendal.
\newblock {\em Stochastic {D}ifferential {E}quations: an introduction with
  applications.}
\newblock Springer Science and Business Media, 2013.

\bibitem{PavliotisStuart:03}
G.~Pavliotis and A.~Stuart.
\newblock White noise limits for inertial particles in a random field.
\newblock {\em Multiscale Model Simul.}, 1:527--553, 2003.

\bibitem{PavliotisStuart:05}
G.~Pavliotis and A.~Stuart.
\newblock Periodic homogenization for inertial particles.
\newblock {\em Physica D}, 204:161--187, 2005.

\bibitem{PavliotisStuart:07}
G.~Pavliotis and A.~Stuart.
\newblock Homogenization for inertial particles in a random flow.
\newblock {\em Commun Math Sci.}, 5:507--531, 2007.

\bibitem{Stuart:08}
G.~Pavliotis and A.~Stuart.
\newblock {\em Multiscale methods: averaging and homogenization}.
\newblock Springer Science and Business Media, 2008.

\bibitem{StuartZygalakis:09}
G.~Pavliotis, A.~Stuart, and K.~Zygalakis.
\newblock Calculating effective diffusivities in the limit of vanishing
  molecular diffusion.
\newblock {\em J. Comput. Phys.}, 228:1030--1055, 2009.

\bibitem{strang:68}
G.~Strang.
\newblock On the construction and comparison of difference schemes.
\newblock {\em SIAM J. Numer. Anal.}, 5:506--517, 1968.

\bibitem{WangXinZhang:18}
Z.~J. Wang, J.~Xin, and Z.~W. Zhang.
\newblock Computing effective diffusivity of chaotic and stochastic flows using
  structure-preserving schemes.
\newblock {\em SIAM Journal on Numerical Analysis}, 56(4):2322--2344, 2018.

\bibitem{xin2016periodic}
J.~Xin, Y.~Yu, and A.~Zlatos.
\newblock Periodic orbits of the abc flow with ${A}={B}={C}=1$.
\newblock {\em SIAM Journal on Mathematical Analysis}, 48(6):4087--4093, 2016.

\bibitem{JackXin:15}
P.~Zu, L.~Chen, and J.~Xin.
\newblock A computational study of residual {K}{P}{P} front speeds in
  time-periodic cellular flows in the small diffusion limit.
\newblock {\em Physica D}, 311:37--44, 2015.

\end{thebibliography}

\end{document}